\documentclass[11pt]{amsart} 
\usepackage[left=3cm,right=3cm,top=3cm,bottom=3cm]{geometry} 
\usepackage{amsmath,amssymb} 
\usepackage{amsmath, amsthm,amssymb,amsfonts,latexsym}
\usepackage{amsthm}
\usepackage{tikz}
\usepackage{lipsum}
\usepackage{color}

\usepackage{hyperref}
\usepackage[capitalise]{cleveref}
\usetikzlibrary{patterns}
\usetikzlibrary{intersections}
\usetikzlibrary{arrows,positioning}
\usepackage{graphicx}
\usepackage{pgfplots}
 \pgfplotsset{compat=1.14}
\usetikzlibrary{plotmarks}
  \usepgfplotslibrary{patchplots}
  \usepackage{grffile}
\usepgfplotslibrary{fillbetween}
\usepackage{float}
\usepackage{caption}
\usepackage{subcaption}

\newcommand*{\rom}[1]{\expandafter\@slowromancap\romannumeral #1@}

\newcommand{\gettikzxy}[3]{%
  \tikz@scan@one@point\pgfutil@firstofone#1\relax
  \edef#2{\the\pgf@x}%
  \edef#3{\the\pgf@y}%
}

\def\namedlabel#1#2{\begingroup
	#2%
	\def\@currentlabel{#2}%
	\phantomsection\label{#1}\endgroup
}

\newcommand{\sB}{\mathcal B}

\newcommand{\sD}{\mathcal D}
\newcommand{\sE}{\mathcal E}
\newcommand{\sF}{\mathcal F}

\newcommand{\sI}{\mathcal I}

\newcommand{\sM}{\mathcal M}

\newcommand{\sP}{\mathcal P}

\newcommand{\sT}{\mathcal T}
\newcommand{\sU}{\mathcal U}

\newcommand{\sX}{\mathcal X}
\newcommand{\sY}{\mathcal Y}

\newcommand{\R}{\mathbb R}
\newcommand{\E}{\mathbb E}

\newcommand{\F}{\mathbb F}

\newcommand{\Prob}{\mathbb P}

\newcommand{\Leb}{\mbox{Leb}}
\newtheorem{thm}{Theorem}[section]
\newtheorem{theorem}{Theorem}[section]

\newtheorem{proposition}{Proposition}[section]

\newtheorem{example}{Example}[section]
\newtheorem{lemma}{Lemma}[section]

\newtheorem{cor}{Corollary}[section]
\newtheorem{remark}{Remark}[section]

\newtheorem{defn}{Definition}[section]

\usepackage[foot]{amsaddr}
\makeatletter
\renewcommand{\email}[2][]{%
	\ifx\emails\@empty\relax\else{\g@addto@macro\emails{,\space}}\fi%
	\@ifnotempty{#1}{\g@addto@macro\emails{\textrm{(#1)}\space}}%
	\g@addto@macro\emails{#2}%
}

\makeatother
\makeatletter
\@namedef{subjclassname@2020}{%
	\textup{2020} Mathematics Subject Classification}
\makeatother
\numberwithin{equation}{section}

\begin{document}
\title[Supermartingale shadow couplings: the decreasing case]{Supermartingale shadow couplings: the decreasing case}

\author{Erhan Bayraktar}
\address{Department of Mathematics, University of Michigan}
\email{erhan@umich.edu}

\author{Shuoqing Deng}
\address{Department of Mathematics, Hong Kong University of Science and Technology}
\email{masdeng@ust.hk}

\author{Dominykas Norgilas}
\email{dnorgila@umich.edu}
\address{Department of Mathematics, University of Michigan}
\thanks{E. Bayraktar is partially supported by the National Science Foundation under grant  DMS-2106556 and by the Susan M. Smith chair.} 



\begin{abstract}
For two measures $\mu$ and $\nu$ that are in convex-decreasing order, Nutz and Stebegg (Canonical supermartingale couplings, Ann. Probab., 46(6):3351--3398, 2018) studied the optimal transport problem with supermartingale constraints and introduced two canonical couplings, namely the increasing and decreasing transport plans, that are optimal for a large class of cost functions. In the present paper we provide an explicit construction of the decreasing coupling $\pi^D$ by establishing a Brenier-type result: (a generalised version of) $\pi^D$ concentrates on the graphs of two functions.

Our construction is based on the concept of the supermartingale \textit{shadow} measure and requires a suitable extension of the results by Juillet (Stability of the shadow projection and the left-curtain coupling, Ann. Inst. H. Poincaré Probab. Statist., 52(4):1823--1843, November 2016) and Beiglb\"ock and Juillet (Shadow couplings, Trans. Amer. Math. Soc., 374:4973--5002, 2021) established in the martingale setting. In particular, we prove the stability of the supermartingale shadow measure with respect to initial and target measures $\mu,\nu$, introduce an infinite family of lifted supermartingale couplings that arise via shadow measure, and show how to explicitly determine the `martingale points' of each such coupling.
\end{abstract}

\maketitle
\tableofcontents
\section{Introduction}\label{sec:intro}
The classical optimal transport (OT) problem is to find a joint law $\pi$ of random variables $X\sim\mu$ and $Y\sim\nu$ that, for a given $c:\R^2\to\R$, minimises the total expected cost $\E^\pi[c(X,Y)]$. Any such joint law corresponds to a measure on $\R^2$, with first and second marginals $\mu$ and $\nu$, respectively, and is called a transport plan from $\mu$ to $\nu$ (or a coupling of $\mu$ and $\nu$). Let $\Pi(\mu,\nu)$ be the set of all such couplings. It is often convenient to express a coupling $\pi\in\Pi(\mu,\nu)$ via its disintegration with respect to the first marginal $\mu$: $\pi(dx,dy)=\mu(dx)\pi_x(dy)$ where $(\pi_x)_{x\in\R}$ is a $\mu$-almost surely unique family of probability kernels.

The cornerstone result in $\R^d$, and with an Euclidean cost $c(x,y)=\lvert x-y\lvert^2$, is Brenier's theorem (see Brenier \cite{Brenier:87} and R\"{u}schendorf and Rachev \cite{RuRach:90}): under some regularity conditions on the initial measure $\mu$, the optimal coupling takes the form $\pi(dx,dy)=\mu(dx)\delta_{\psi(x)}(dy)$, where $\psi:=\nabla\phi$ is the gradient of a convex function $\phi$. In dimension one, the supporting function is non-decreasing and the optimal coupling coincides with the Hoeffding-Fr\'echet (or quantile) coupling $\pi^{HF}$, which, in the case $\mu$ is continuous, is given by $\pi^{HF}(dx,dy)=\mu(dx)\delta_{G_\nu(F_\mu(x))}(dy)$, where $G_\nu$ and $F_\mu$ are the quantile function of $\nu$ and the cumulative distribution function of $\mu$, respectively. An important feature of $\pi^{HF}$ is that it is optimal for a large class of cost functions (essentially those that satisfy the Spence-Mirrlees condition $c_{xy}>0$).

In the last decade, there has been a significant interest in the OT problems where the coupling $\pi$ is required to constitute a martingale. In particular, in the martingale optimal transport (MOT) one still seeks to minimise (or maximise) the total expected cost (or payoff) $\int_{\R^2}cd\pi$, but only over the set of martingale transport plans: $\pi\in\Pi(\mu,\nu)$ is a martingale coupling, and we write $\pi\in\Pi_M(\mu,\nu)\subseteq \Pi(\mu,\nu)$, if $
\int_\R y\pi_x(dy)= x$ for $\mu$-a.e. $x\in\R$ (or equivalently, $\pi\in\Pi_M(\mu,\nu)$ if $\E^\pi[Y\lvert X]=X$). Such problems arise naturally in the context of model-independent mathematical finance (see Beigl\"ock et al. \cite{BeiglbockHenryLaborderePenkner:13} and Galichon et al. \cite{GalichonHenryLabordereTouzi:14}) and have important consequences for martingale inequalities (see Beiglb\"{o}ck and Nutz \cite{BeiglbockNutz:14}, Henry-Labord{\`e}re et al. \cite{HLabordereOblojSpoida:16}, Ob{\l}{\'o}j et al. \cite{OblojSpoidaTouzi:15}) and the Skorokhod embedding problem (see Beiglb\"{o}ck et al. \cite{BeiglbockCoxHuesmann:17}, K\"{a}llblad et al. \cite{KallbladTanTouzi:17}).

The first explicit solutions to the MOT problem is due to Hobson and Neuberger \cite{HobsonNeuberger:12} and Hobson and Klimmek \cite{HobsonKlimmek:15} where the authors showed how construct couplings $\pi^{HN}$ and $\pi^{HK}$ that maximise and minimise $\E^\pi[\lvert Y-X\lvert]$, receptively. It is not known, however, whether these couplings are optimal for non-Euclidean costs functions. A more general result was obtained by Beiglb\"ock et al. \cite{BeiglbockJuillet:16}. Using an extension of the cyclical monotonicity from the classical OT setting, the authors introduced the left-monotone martingale coupling and baptised it as the \textit{left-curtain} transport plan $\pi^{lc}$. It was shown that such left-monotone coupling exists, is unique and also optimal for a particular class of cost functions. Henry-Labord\'ere and Touzi \cite{HenryLabordereTouzi:16} extended the results of \cite{BeiglbockJuillet:16} to show that $\pi^{lc}$ is optimal for even wider class of payoffs (namely those that satisfy the martingale Spence-Mirrlees condition $c_{xyy}<0$), and in the case when the initial measure $\mu$ is atom-less, provided an explicit construction using (coupled) ordinary differential equations. A general construction for arbitrary $\mu$ and $\nu$ was recently obtained by Hobson and Norgilas \cite{HobsonNorgilas:21}. Several other properties and extensions of $\pi^{lc}$ where further investigated in Beiglb\"{o}ck et al. \cite{{BeiglbockHenryLabordereTouzi:17},{BeiglbockCox:17},BHN:20}, Juillet \cite{{Juillet:16},Juillet:18}, Hobson and Norgilas \cite{HobsonNorgilas:17}, Nutz et al. \cite{{NutzStebegg:18},NutzStebeggTan:17}, Campi et al. \cite{Campi:17}, Henry-Labord\`{e}re et al.~\cite{HenryLabordereTanTouzi:16} and Br\"{u}ckerhoff at al.~\cite{BHJ:20}.

In this paper we shall study the couplings that constitute a supermartingale: $\Pi_S(\mu,\nu)\subseteq\Pi(\mu,\nu)$ is the set of all supermartingale couplings, and we write $\pi\in\Pi_S(\mu,\nu)$ if 
$$
\int_\R y\pi_x(dy)\leq x,\quad\textrm{for }\mu\textrm{-a.e. }x\in\R.
$$
By a classical result of Strassen \cite{Strassen:65}, $\Pi_S(\mu,\nu)$ is non-empty if and only if $\mu$ and $\nu$ are in convex-decreasing order (denoted by $\mu\leq_{cd}\nu$), i.e., $\int_\R fd\mu\leq\int_\R fd\nu$ for all convex and non-increasing $f:\R\to\R$ (if the inequality holds for all convex functions, then $\mu$ and $\nu$ are in convex order, i.e., $\mu\leq_c\nu$, and $\Pi_M(\mu,\nu)\neq\emptyset$). Note that $\Pi_M(\mu,\nu)\subseteq\Pi_S(\mu,\nu)$. In particular, if $\mu\leq_{cd}\nu$ and with equal means then $\Pi_M(\mu,\nu)=\Pi_S(\mu,\nu)$. However, even when $\mu\leq_{cd}\nu$ holds, it is a non-trivial problem to construct particular supermartingale transport plans.

Similarly as in the contexts of OT and MOT, for a given $c:\R^2\to\R$ it is natural to seek for elements $\pi\in\Pi_S(\mu,\nu)$ that minimise/maximise the total expected cost $\int_{\R^2}c(x,y)\pi(dx,dy)$. In this setting Nutz and Stebegg \cite{NutzStebegg:18} introduced two canonical supermartingale couplings, namely the increasing coupling $\pi^I$ and decreasing coupling $\pi^D$. Our main interest in this article is the latter. Nutz and Stebegg \cite{NutzStebegg:18} proved that $\pi^D$ is a unique supermartingale coupling that can be equivalently characterised by any, and then all of the following properties (see Definition \ref{defn:decreasing_coupling}): the optimality (in terms of minimization problem) for a large class of cost functions (essentially those that satisfy $c_{xy}<0$ and $c_{xyy}>0$), the monotonicity of the support, and the (stochastic) order-theoretic minimality. (A similar characterisation is valid for $\pi^I$ as well.) While this can be viewed as an existence result, our aim is to provide an explicit construction of $\pi^D$. In particular, and in the light of Brenier's theorem in OT, our goal is to construct functions on the graph of which the decreasing coupling $\pi^D$ concentrates. (Indeed, all of the aforementioned martingale couplings and the increasing supermartingale coupling $\pi^I$ are constructed in such way.)

The main ingredient in our construction is the so-called \textit{shadow} measure, introduced by Beiglb\"ock and Juillet \cite{BeiglbockJuillet:16} in the martingale setting and later extended by Nutz and Stebegg \cite{NutzStebegg:18} to supermartingales. For $\eta\leq\mu$ the (supermartingale) shadow of $\eta$ in $\nu$, denoted by $S^\nu(\eta)$, is the measure satisfying $\eta\leq_{cd} S^\nu(\eta)\leq\nu$ and $S^\nu(\eta)\leq_{cd}\chi$ for all $\eta\leq_{cd}\chi\leq\nu$ (if one replaces $\leq_{cd}$ by $\leq_c$ then the notion of the martingale shadow measure $S_M^\nu(\eta)$ is recovered). In particular, $S^\nu(\eta)$ is the smallest (with respect to $\leq_{cd}$) measure within $\nu$ to which $\eta$ can be mapped to in a supermartingale way. Our interest in $S^\nu(\eta)$ lies in the fact that the decreasing supermartingale coupling $\pi^D$ can be defined as a unique coupling such that, for each $x\in\R$, $\pi^D\lvert_{[x,\infty)\times\R}$ has
first marginal $\mu\lvert_{[x,\infty)}$ and second marginal $S^\nu(\mu\lvert_{[x,\infty)})$. On the other hand, Bayraktar et al. \cite{BDN:21} (among other things) showed how to construct a potential function of the shadow measure $S^\nu(\eta)$. (A potential function is a convex function whose second derivative uniquely identifies the underlying measure.) One of our main contributions is to show how, given the graph of the potential function of $S^\nu(\eta)$, to identify a pair of candidate functions that support the decreasing supermartingale coupling $\pi^D$.

An idea of constructing canonical couplings via corresponding potential functions is not new and dates back to Hobson and Klimmek \cite{HobsonKlimmek:15} at least (such approach is often taken to obtain particular solutions to the Skorokhod embedding problem; see Ob\l{}\'oj \cite{Obloj:survey} and Hobson \cite{Hobson:survey} for an overview). On the other hand, Hobson and Norgilas \cite{HobsonNorgilas:21} were first to show that the functions that support the left-curtain coupling $\pi^{lc}$ can be identified from the potential function of the (martingale) shadow measure ${S}_M^\nu(\mu\lvert_{(-\infty,x]})$, $x\in\R$. Similarly, one can construct the \textit{right-curtain} coupling $\pi^{rc}$ (the symmetric counterpart of $\pi^{lc}$) by studying the potential functions of ${S}_M^\nu(\mu\lvert_{[x,\infty)})$, $x\in\R$. In the supermartingale setting, Bayraktar et al. \cite{BDN:21} used this approach and constructed a pair of functions that support the increasing supermartingale coupling $\pi^I$.

The increasing supermartingale transport plan $\pi^I$ is obtained by working from left to right (i.e., by mapping, for each $x\in\R$, $\mu\lvert_{(-\infty,x]}$ to $S^\nu(\mu\lvert_{(-\infty,x]})$) and thus can be viewed as a supermartingale counterpart of $\pi^{lc}$. (Indeed if $\mu\leq_{cd}\nu$ and with equal means then $\mu\leq_c\nu$ and $S^\nu(\eta)=S_M^\nu(\eta)$ for all $\eta\leq\mu$.) One of the main achievements of Bayraktar et al. \cite{BDN:21} is that the authors showed how to obtain $x^*\in\R$, such that (under $\pi^I$) $\mu\lvert_{(-\infty,x^*]}$ is embedded in $\nu$ via the martingale shadow measure $S_M^\nu(\mu\lvert_{(-\infty,x^*]})$, while the remaining mass $(\mu-\mu\lvert_{(-\infty,x^*]})$ is mapped to $\nu$ via supermartingale shadow $S^\nu(\cdot)$. In particular, $\pi^I\lvert_{(-\infty,x^*]\times\R}=\pi^{lc}\lvert_{(-\infty,x^*]\times\R}$ while $\pi^I\lvert_{(x,\infty)\times\R}$ corresponds to the (deterministic) antitone coupling $\pi^{AT}$ (a symmetric counterpart of the quantile coupling $\pi^{HF}$).

On the other hand, the decreasing coupling $\pi^D$ is a supermartingale variant of the right-curtain coupling $\pi^{rc}$. Indeed, both transport plans are constructed by embedding $\mu\lvert_{[x,\infty)}$ in $\nu$ via $S^\nu(\cdot)$ and $S^\nu_M(\cdot)$, respectively. More precisely, and as already observed in Nutz and Stebegg \cite{NutzStebegg:18}, $\pi^D$ corresponds to $\pi^{rc}$ on a specific part of the state space (the so-called `martingale points') and is equal to the quantile coupling $\pi^{HF}$ elsewhere. Even though we can explicitly determine the `martingale points' of $\pi^D$ (see Section \ref{sec:martingale_points}), there could be a countably many regime switching points where $\pi^D$ alternates between $\pi^{rc}$ and $\pi^{HF}$ (compare this with a unique regime switching point in the case of $\pi^I$). Therefore (due to the lack of symmetry between $\pi^I$ and $\pi^D$) the pasting arguments of Bayraktar et al. \cite{BDN:21} seem to be hard to adapt. Instead, we relate our construction to the notion of the lifted shadow couplings introduced by Beiglb\"ock and Juillet \cite{BeiglbockJuillet:16s} in the martingale setting.

In the case of lifted couplings the idea is to consider a \textit{lift} of $\mu$, given by $\hat\mu\in\Pi(\lambda,\mu)$, where $\lambda$ is the Lebesgue measure on $[0,1]$, and then to construct a lifted (martingale) transport plan $\hat\pi$ that has first and second marginals $\hat\mu$ and $\nu$, respectively. Note that each such $\hat\pi$ is a measure on $[0,1]\times\R\times\R$. If we disintegrate $\hat\pi$ with respect to $\hat\mu$, so that $\hat\pi(du,dx,dy)=\hat\mu(du,dx)\hat\pi_{(u,x)}(dy)$, then the martingale condition reads $\int_\R y\hat\pi_{(u,x)}(dy)=x$ for $\hat\mu$-a.e. $(u,x)\in[0,1]\times\R$. A corresponding martingale coupling on $\R^2$ is then recovered by integrating out the Lebesgue measure. Beiglb\"ock and Juillet \cite{BeiglbockJuillet:16s} showed that for each lift $\hat\mu$ there exists a unique lifted martingale \textit{shadow coupling} $\hat\pi$ that, for each $u\in[0,1]$, maps $\hat\mu\lvert_{[0,u]\times\R\times\R}$ to $S_M^\nu(\hat\mu\lvert_{[0,u]\times\R\times\R})$. (Then in order to obtain a (a version of a) particular martingale coupling, e.g., $\pi^{lc}$ or $\pi^{rc}$, one just needs to choose an appropriate lift $\hat\mu$.) More precisely, represent $\hat\pi$ as $\hat\pi(du,dx,dy)=du\hat\mu_{u,\cdot}(dx)\hat\pi_{(u,x)}(dy)$, where $(\hat\mu_{u,\cdot})_{u\in[0,1]}$ corresponds to the disintegration of $\hat\mu$ with respect to $\lambda$, while $(\hat\pi_{(u,x)})_{(u,x)\in[0,1]\times\R}$ is the disintegration of $\hat\pi$ with respect to $\hat\mu$. Then one of the main insights of Beiglb\"ock and Juillet \cite{BeiglbockJuillet:16s} is that in fact the kernel $\hat\pi_{(u,x)}$ corresponds to the hitting coupling (of $\hat\mu_{u,\cdot}$ to a suitably defined subset of $\R$) introduced by Kellerer \cite{kellerer73}. A crucial observation for our purposes is that when $\hat\mu_{u,\cdot}$ is a point mass, then the hitting coupling $\hat\pi_{(u,x)}$ is a measure concentrated on at most two points. Our aim is to generalise these results in the supermartingale setting, and then show how, using the potential function of the shadow measure $S^\nu(\cdot)$, to identify the points on which the \textit{supermartingale hitting coupling} concentrates. Consequently, this will allow us to recover the supporting functions (and thus also an explicit construction) of $\pi^D$.

There are several main results in this paper, for each of which there is a dedicated section. First, in Section \ref{sec:shadow} we prove the stability of the supermartingale shadow measure $S^\nu(\eta)$ with respect to the given data $\eta\leq\mu\leq_{cd}\nu$; see Theorem \ref{thm:shadowW} (this generalises the results of Juillet \cite{Juillet:16} obtained for martingales). This is a crucial ingredient in the proof of existence and uniqueness result (see Theorem \ref{thm:shadowCouplings}) regarding the lifted supermartingale shadow couplings (this generalises the corresponding result obtained by Beiglb\"ock and Juillet \cite{BeiglbockJuillet:16s} in the martingale setting); this is presented in Section \ref{sec:mainExistence}. In Section \ref{sec:martingale_points} we show how to explicitly determine the `martingale points' of each lifted shadow coupling of Theorem \ref{thm:shadowW}; see Proposition \ref{prop:martingale_points}. Finally, in Section \ref{sec:geometric} we link the potential function of the shadow measure and the transition kernel of a lifted shadow coupling, and show how to explicitly construct the increasing supermartingale coupling $\pi^D$; see Theorem \ref{thm:construction}.

\section{Preliminaries}\label{sec:prelims}
\subsection{Spaces of measures and related notions}\label{sec:measureSpaces}
For $d\geq 1$, we denote by $\sM^d$ (resp. $\sP^d$) the space of positive measures (resp. probability measures) on $\R^d$ with finite first moments. In the case $d=1$ we write $\sP=\sP^1$ and $\sM=\sM^1$.

The support of a measure $\eta\in\sM^d$ is denoted by $supp(\eta)$. It is the smallest closed (Borel) set $E\subseteq \R^d$ with $\eta(\R^d\setminus E)=0$.

We use $\lambda_I$ to denote the restriction of the Lebesgue measure (on $\R$) to an interval $I\subseteq\R$. In the case $I=[0,1]$, we write $\lambda=\lambda_{[0,1]}$.

 Given a measure $\eta\in\sM$ (not necessarily a probability measure), define $\bar{\eta} = \int_\R x \eta(dx)$ to be the first moment of $\eta$ (and then $\bar\eta/\eta(\R)$ is the barycentre of $\eta$). Let $\sI_\eta$ be the smallest interval containing the support of $\eta$, and let $\{ \ell_\eta, r_\eta \}$ be the endpoints of $\sI_\eta$. If $\eta$ has an atom at $\ell_\eta$ then $\ell_\eta$ is included in $\sI_\eta$, and otherwise it is excluded, and similarly for $r_\eta$.

For $\mu\in\sM$, the right-continuous cumulative distribution function $F_\mu:\R\to[0,\mu(\R)]$ is defined by $F_\mu(x):=\mu((-\infty,x])$, $x\in\R$. A quantile function of $\mu$, i.e., a generalised inverse of $F_\eta$, is denoted by $G_\eta:[0,\eta(\R)]\mapsto\R$. There are two canonical versions of $G_\eta$: the left-continuous and right-continuous versions correspond to $G_\eta^-(u)=\sup\{k\in\R: F_\eta(k)<u\}$ and $G^+_\eta(u)=\inf\{k\in\R: F_\eta(k) > u\}$, for $u\in[0,\eta(\R)]$, respectively. However any $G$ with $G_\eta^-(u)\leq G(u)\leq G^+_\eta(u)$, for all $u\in[0,\eta({\R}]$, is still called a quantile function of $\eta$, which is motivated by the fact that for any such $G$ we have that Law$(G(U))=\eta/\eta(\R)$, where $U\sim U[0,\eta(\R)]$. (Note that $G_\eta$ may take values $-\infty$ and $\infty$ at the left and right end-points of $[0,\eta(\R)]$, respectively.) 

\subsubsection{Potential functions} For $\alpha \geq 0$ and $\beta \in \R$ let $\sD^\uparrow(\alpha, \beta)$ denote the set of non-negative, non-decreasing and convex functions $f:\R \mapsto \R_+$ such that
\[ \lim_{ z \downarrow -\infty}   f(z)  =  0, \hspace{10mm} \lim_{z \uparrow \infty} \{ f(z) - (\alpha z- \beta) \}   =0. \]
Then, when $\alpha = 0$, $\sD^\uparrow(0,\beta)$ is empty unless $\beta = 0$ and then $\sD^\uparrow(0,0)$
contains one element, the zero function. Similarly, let $\sD^\downarrow(\alpha, \beta)$ denote the set of non-negative, non-increasing and convex functions $f:\R \mapsto \R_+$ such that
\[ \lim_{ z \downarrow -\infty}  \{ f(z)  -(\beta-\alpha z) \} =  0, \hspace{10mm} \lim_{z \uparrow \infty}  f(z)    =0. \]

For $\eta\in\sM$, define the functions $P_\eta,C_\eta : \R \mapsto \R_+$ by
\begin{equation*}
P_\eta(k) := \int_{\R} (k-x)^+ \eta(dx),\quad k\in\R,\hspace{10mm}C_\eta(k) := \int_{\R} (x-k)^+ \eta(dx),\quad k\in\R,
\end{equation*}
respectively. Then $P_\eta(k) \geq 0 \vee  (\eta(\R) k - \overline{\eta} )$ and $C_\eta(k) \geq 0 \vee (\overline{\eta} - \eta(\R)k)$. Also, the Put-Call parity holds: $C_\eta(k) - P_\eta(k) = (\overline{\eta}- \eta(\R)k)$, $k\in\R$.

The following properties of $P_\eta$ (resp. $C_\eta$) can be found in Chacon~\cite{Chacon:77}, and Chacon and Walsh~\cite{ChaconWalsh:76}: $P_\eta \in \sD^\uparrow(\eta(\R), \overline{\eta})$ (resp. $C_\eta \in \sD^\downarrow(\eta(\R), \overline{\eta})$) and $\{k : P_{\eta}(k) > (\eta(\R)k - \overline{\eta})^+  \} =  \{k : C_{\eta}(k) > (\overline{\eta}- \eta(\R)k)^+\} =(\ell_\eta,r_\eta)$. Conversely (see, for example, Proposition 2.1 in Hirsch et al. \cite{peacock}),  if $h \in \sD^\uparrow(k_m,k_f)$ for some numbers $k_m \geq 0$ and $k_f\in\R$ (with $k_f = 0$ if $k_m=0$), then there exists a unique measure $\eta\in\sM$, with total mass $\eta(\R)=k_m$ and mean $\overline{\eta}=k_f$, such that $h=P_{\eta}$. In particular, $\eta$ is uniquely identified by the second derivative of $h$ in the sense of distributions. Furthermore, $P_\eta$ and $C_\eta$ are related to the potential $U_\eta$, defined by
\begin{equation*}
U_\eta(k) : =  - \int_{\R} |k-x| \eta(dx),\quad k\in\R,
\end{equation*}
by $-U_\eta=C_\eta+P_\eta$. We will call $P_\eta$ (and $C_\eta$) a modified potential. Finally note that all three second derivatives $C^{\prime\prime}_{\eta},P^{\prime\prime}_{\eta}$ and $-U_\eta^{\prime\prime}/2$ identify the same underlying measure  $\eta$.

\subsubsection{Wasserstein distance}\label{sec:Wasserstein}
For $(\mu_n)_{n\geq1},\mu \in\sM$ with $\mu_n(\R)=\mu(\R)$ for all $n\geq1$, we use the notation $\mu_n\xrightarrow{w}\mu$ for the usual weak convergence of measures, i.e., $\mu_n\xrightarrow{w}\mu$ if $\int_\R fd\mu_n\to\int_\R fd\mu$ for all continuous and bounded $f:\R\to\R$.

For $\mu,\nu\in\sP$, the Wasserstein-1 distance is defined by
\begin{equation}\label{eq:Wass}
W(\mu,\nu)=\sup_{f\in Lip(1)}\left\lvert \int fd\mu -\int fd\nu\right\lvert,
\end{equation}
where the supremum is taken over all 1-Lipschitz functions $f:\R\to\R$. It endows $(\sP,W)$ with $\sT_1$, the usual topology for probability measures with finite first moments (a sequence of measures $(\mu_n)_{n\geq1}$ converges to $\mu$ w.r.t. $\sT_1$, and we write $\mu_n\to\mu$, if $\int_\R fd\mu_n\to\int_\R fd\mu$ for every continuous $f:\R\to\R$ with at most linear growth). Moreover, if $\mu_n\to\mu$ (or, equivalently, $W(\mu_n,\mu)\to0$) then $\mu_n\xrightarrow{w}\mu$, while the converse is true if the first moments also converge (see, for example, Villani \cite[Theorem 6.9]{villani:08}). For a fixed $m\in\R$, we will often work with a subspace $\{\eta\in\sM:\eta(\R)=m\}\subseteq\sM$, and in this case we still consider the distance $W$ with the same definition.

Note that, for each $t\in\R$, $f_t:x\in\R\to-\lvert x-t\lvert\in\R$ belongs to $Lip(1)$, and therefore, if $\mu_n\to\mu$ in $\sM$, then the potential functions $U_{\mu_n} $ converge to $U_\mu$ pointwise. The reverse implication does not hold in general and we need additional assumptions on the given data (for example, it is enough for all the measures to have the same mass and mean, see Hirsch and Roynette \cite[Proposition 2.3]{peacock}).

By Kantorovich duality theorem (see  Villani \cite{villani:08}), and in the case $\mu,\nu\in\sP$, one has an alternative definition of the Wasserstein-1 distance given by
$$
W(\mu,\nu)=\inf_{\pi}\int_{\R^2}\lvert y-x\lvert d\pi(x,y),
$$
where the infimum is over all $\pi\in\sP^2$ with marginals $\mu$ and $\nu$. In particular, the infimum is attained by the Hoeffding-Fr\'echet (or \textit{quantile}) coupling $\pi^{HF}\in\sP^2$, defined by
$$
\pi^{HF}(A\times B)=\lambda(\{u\in[0,1]:G_\mu(u)\in A,G_\nu(u)\in B\}),\quad A,B\in\sB(\R).
$$
Then
\begin{equation}\label{eq:WassQuantile}
W(\mu,\nu)=\int^1_0\lvert G_\nu(u)-G_\mu(u)\lvert du,
\end{equation}
which, in the case when $\mu,\nu\in\sM$ are with equal mass, generalises to $$W(\mu,\nu)=\int^{\mu(\R)=\nu(\R)}_0\lvert G_\nu(u)-G_\mu(u)\lvert du.$$

\subsection{Stochastic orders and supermartingale couplings}\label{sec:stochasticOrder}
For $\eta,\chi\in\sM$, we write $\eta\leq\chi$ if $\eta(A) \leq \chi(A)$ for all Borel measurable subsets $A$ of $\R$, or equivalently if
\begin{equation*}
\int fd\eta\leq\int fd\chi,\quad \textrm{for all non-negative }f:\R\mapsto\R.
\end{equation*}
Since $\eta$ and $\chi$ can be identified as second derivatives of the potential functions $P_\eta$ and $P_\chi$ , we have $\eta\leq\chi$ if and only if $P_\chi-P_\eta$ is convex, i.e., $P_\eta$ has a smaller curvature than $P_\chi$.

Two measures $\eta,\chi\in\sM$ with equal mass are in stochastic order, and we write $\eta\leq_{sto}\chi$, if 
$$
\int fd\eta\leq\int fd\chi,\quad\textrm{for all non-decreasing }f:\R\mapsto\R.
$$
In particular, $\eta\leq_{sto}\chi$ if and only if $F_\eta\geq F_\chi$ on $\R$, or equivalently, $G_\eta\leq G_\chi$ on $[0,1]$.

Two measures $\eta,\chi\in\sM$ are in convex (resp. convex-decreasing) order, and we write $\eta\leq_c\chi$ (resp. $\eta \leq_{cd} \chi$), if
\begin{equation}\label{eq:cd}
\int fd\eta\leq\int fd\chi,\quad\textrm{for all convex (resp. convex and non-increasing) }f:\R\mapsto\R.
\end{equation}
Since we can apply \eqref{eq:cd} to all constant functions, including $f\equiv -1$ and $f\equiv1$, we have that if $\eta\leq_{c}\chi$ (or $\eta\leq_{cd}\chi$) then $\eta(\R)=\chi(\R)$. On the other hand, applying \eqref{eq:cd} to $f(x)=-x$ gives that $\bar{\eta}\geq\bar{\chi}$ whenever $\eta\leq_{cd}\chi$. However, a reversed inequality holds only in the case $\eta\leq_c\chi$ (since $f(x)=x$ is strictly increasing).

Given two probability measures $\eta,\chi$ on Polish spaces $\sX,\sY$, respectively, let ${\Pi}(\eta,\chi)$ be the set of probability measures on $\sX\times\sY$ with the first marginal $\eta$ and second marginal $\chi$.

For $\eta,\chi\in\sP$ let ${\Pi}_S(\eta,\chi)$ be the set of supermartingale couplings of $\eta$ and $\chi$.
Then
\begin{equation*}
{\Pi}_S(\eta,\chi) = \big\{ \pi \in {\Pi}(\eta,\chi) : \mbox{\eqref{eq:martingalepi} holds} \big\},
\end{equation*}
where \eqref{eq:martingalepi} is the supermartingale condition
\begin{equation}
\int_{x \in B} \int_{y \in \R} y \pi(dx,dy) \leq \int_{x \in B} \int_{y \in \R} x \pi(dx,dy) = \int_B x \eta(dx), \quad
\mbox{$\forall$ Borel $B \subseteq \R$}.
\label{eq:martingalepi}
\end{equation}
Equivalently, $\Pi_S(\eta,\chi)$ consists of all transport plans $\pi$ (i.e., elements of ${\Pi}(\eta,\chi)$) such that the disintegration in probability measures $(\pi_x)_{x\in\R}$ with respect to $\eta$ satisfies $\int_\R y\pi_x(dy)\leq x$ for $\eta$-almost every $x$.

The following is classical (see, for example, F\"{o}llmer and Schied \cite[Theorem 2.58]{FS:16}).
\begin{lemma}\label{lem:put_cd}
	Let $\eta,\chi\in\sP$. The following are equivalent:\begin{enumerate}
		\item $\eta\leq_{cd}\chi$,
		\item $\eta(\R)=\chi(\R)$ and $P_\eta\leq P_\chi$ on $\R$,
		\item $\Pi_S(\eta,\chi)\neq\emptyset$.
	\end{enumerate}
\end{lemma}

	If $\eta,\chi\in\sP$ with $\eta\leq_{cd}\chi$, but $\bar\eta=\bar\chi$, then $\Pi_S(\eta,\chi)$ reduces to the set of martingale couplings, denoted by $\Pi_M(\eta,\chi)$ (i.e., elements of $\Pi(\eta,\chi)$ for which \eqref{eq:martingalepi} holds with equality). Indeed, any supermartingale with constant mean is a martingale. In this case $\eta\leq_c\chi$ (see Strassen~\cite{Strassen:65}).

For our purposes in the sequel we need a generalisation of the convex (resp. convex-decreasing) order of two measures. We say $\eta,\chi\in\sM$ are in a \textit{positive} convex (resp. \textit{positive} convex-decreasing) order, and write $\eta\leq_{pc}\chi$ (resp. $\eta\leq_{pcd}\chi$), if $\int fd\eta\leq\int fd\chi$, for all non-negative and convex (resp. non-negative, convex and non-increasing) $f:\R\mapsto\R_+$.
If $\eta\leq_{c}\chi$ (resp. $\eta\leq_{cd}\chi$) then also $\eta\leq_{pc}\chi$ (resp. $\eta\leq_{pcd}\chi$), since non-negative and convex (resp. non-negative, convex and non-increasing) functions are convex (resp. convex and non-increasing). If $\eta\leq\chi$ then both,  $\eta\leq_{pc}\chi$ and $\eta\leq_{pcd}\chi$, since non-negative, convex and non-increasing functions are non-negative and convex, and thus also non-negative. Note that, if $\eta\leq_{pc}\chi$ or $\eta\leq_{pcd}\chi$, then $\eta(\R)\leq\chi(\R)$ (apply the function $f(x)\equiv1$ in the definition of $\leq_{pc}$ and $\leq_{pcd}$). On the other hand, if $\eta(\R)=\chi(\R)$, then $\eta\leq_{pc}\chi$ (resp. $\eta\leq_{pcd}\chi$) is equivalent to $\eta\leq_{c}\chi$ (resp. $\eta\leq_{cd}\chi$).
\begin{example}\label{eg1}
	Let $\eta,\chi\in\sM$ with $\eta\leq_{cd}\chi$ (resp. $\eta\leq_{c}\chi$). Fix a Borel set $B\subseteq\R$, and let $\eta\lvert_B\in\sM$ be a restriction of $\eta$ to $B$. Then $\eta\lvert_B\leq_{pcd}\chi$ (resp. $\eta\lvert_B\leq_{pc}\chi$).
\end{example}
Let $\eta,\chi\in\sM$ with $\eta\leq_{pcd}\chi$, and introduce the set
\begin{equation}\label{eq:M_nu_mu}
\sM^\chi_\eta:=\{\theta\in\sM:\eta\leq_{cd}\theta\leq\chi\}.
\end{equation}
Then $\sM^\chi_\eta$ is the set of target measures of a supermartingale that transports (or embeds) $\eta$ to $\chi$. It is not hard to show that $\sM^\chi_\eta\neq\emptyset$. Indeed, the `left-most' measure $\theta\leq\chi$ of mass $\eta(\R)$, denoted by $\underline\theta=\underline{\theta}_{\eta,\chi}\in\sM$, and defined by
\begin{equation}\label{eq:left-most}
\underline{\theta}=\chi\lvert_{(-\infty,G_\chi(\eta(\R))}+(\eta(\R)-\chi\lvert_{(-\infty,G_\chi(\eta(\R))}(\R))\delta_{G_\chi(\eta(\R))},
\end{equation}
is the largest measure in $\sM^\chi_\eta$ with respect to $\leq_{cd}$ (see Bayraktar et al. \cite[Proposition 3.1]{BDN:21}). Note that $\underline\theta$ does not depend on the choice of the quantile function $G_\chi$. 
\subsubsection{Irreducible decomposition of $\pi\in\Pi_S(\mu,\nu)$}\label{sec:irreducible}
For a pair of measures $\eta,\chi\in\sM$, let the function $D = D_{\eta,\chi}: \R \mapsto \R$ be defined by $D_{\eta,\chi}(k) = P_\chi(k) - P_\eta(k)$. Note that if $\eta,\chi$ have equal mass then $\eta \leq_{cd} \chi$ is equivalent to $D \geq 0$ on $\R$. In particular, $\lim_{k\to-\infty}D(k)=0$ and $\lim_{k\to\infty}D(k)=\overline\eta-\overline\chi\geq0$. Let $\ell_D:=\inf\{k\in\R:D_{\eta,\chi}(k)>0\}$ and $r_D:=\sup\{k\in\R:D_{\eta,\chi}(k)\neq\overline\eta-\overline\chi\}$. (Note that if $\eta\leq_{cd}\chi$ and $\overline\eta=\overline\chi$ then $r_D=\sup\{k\in\R:D_{\eta,\chi}(k)>0\}$.) Let $\sI_D$ be the open interval $(\ell_D,r_D)$ together with $\{\ell_D\}$ if $\ell_D>-\infty$ and $D'(\ell_D+):=\lim_{k\downarrow\ell_D}(D(k)-D(\ell_D))/(k-\ell_D)>0$ and $\{ r_D \}$ if $r_D<\infty$ and $D'(r_D-):=\lim_{k\uparrow\ell_D}(D(k)-D(\ell_D))/(k-\ell_D) \neq0$. Note that, if $\eta\leq_c\chi$ (or equivalently, $\eta\leq_{cd}\chi$ and $\overline\eta=\overline\chi$), then $\ell_\nu\leq\ell_\mu\leq r_\mu\leq r_\nu$ and $\sI_D\subseteq[\ell_\nu,r_\nu]$. On the other hand, if $\eta\leq_{cd}\chi$ then in general we only have that $\ell_\nu\leq\ell_\mu$ and $\ell_\nu\leq\ell_D$.

It is well know (see Hobson~\cite[page 254]{Hobson:98b} or Beiglb\"{o}ck and Juillet ~\cite[Section A.1]{BeiglbockJuillet:16}) that, if $D_{\eta,\chi}(x)=0$ for some $x$, then in any martingale coupling of $\eta$ and $\chi$ no mass can cross $x$. More precisely, if $\eta \leq_{c} \chi$, $\pi \in {\Pi_M}(\eta,\chi)$ and $x$ is such that $D(x)=0$, then we have $\pi((-\infty,x),(x,\infty)) + \pi((x,\infty),(-\infty,x))=0$.

In the supermartingale case with $\eta\leq_{cd}\chi$, define $x^*:=x^*_{\eta,\chi}:=\sup\{x\in\R: D(x)=0\}$. Then in any supermartingale coupling $\pi\in\Pi_S(\eta,\chi)$, if $D(x)=0$ for some $x\leq x^*$  then no mass can cross $x$ and $\pi$ is a martingale on $(-\infty,x]\times\R$, and $x^*$ is the largest such $x$ (see Nutz and Stebegg \cite{NutzStebegg:18}). In particular we can separate the problem of constructing supermartingale couplings of $\eta$ to $\chi$ into a pair of subproblems involving mass to the left and right of $x^*$, respectively, always taking care to allocate mass of $\chi$ at $x^*$ appropriately. More generally, if there are multiple $\{ x_j \}_{j\geq1}$ (with $x_j\leq x^*$ and) with $D_{\eta,\chi}(x_j)=0$, then we can divide the problem into a sequence of `irreducible' problems,  each taking place on an interval $\sI_i$ such that $D>0$ on the interior of $\sI_i$ and $D=0$ at the endpoints. All mass starting in a given interval is transported to a point in the same interval. Moreover, by the martingale property, any mass starting at a finite endpoint of $\sI_i$ (that is smaller than $x^*$) must stay there. Putting this together we may restrict attention to intervals $I$ on which $D>0$ (with $\lim_{x \rightarrow e_I} D(x)=0$ at endpoints $e_I$ of $I$), and we may assume that the starting law has support within the interior of $I$ and the target law has support within the closure of $I$ (and $I$ is the smallest set with this last property). This is summarized in the following result.
\begin{lemma}[{Nutz and Stebegg \cite[Proposition 3.4]{NutzStebegg:18}}]\label{lem:irreducible}
Let $\mu,\nu\in\sP$ with $\mu\leq_{cd}\nu$. Define $I_0:=(x^*,+\infty)$, let $(I_k)_{k\geq1}$ be the open components of $\{D>0\}\cap(-\infty,x^*)$ and set $I_{-1}:=\R\setminus\bigcup_{k\geq0}I_k$. Let $\mu_k:=\mu\lvert_{I_k}$ for $k\geq-1$, so that $\mu=\sum_{k\geq-1}\mu_k$.

Then there exists a unique decomposition $\nu=\sum_{k\geq-1}\nu_k$ such that
$$
\mu_{-1}=\nu_{-1},\quad \mu_0\leq_{cd}\nu_0\quad\textrm{and}\quad\mu_k\leq_c\nu_k\quad\textrm{for all }k\geq1.
$$
Furthermore, any $\pi\in\Pi_S(\mu,\nu)$ admits a unique decomposition $\pi=\sum_{k\geq-1}\pi_k$ such that $\pi_0\in\Pi_S(\mu_0,\nu_0)$ and $\pi_k\in\Pi_M(\mu_k,\nu_k)$ for all $k\neq0$.
\end{lemma}

\subsection{Lifted couplings}\label{sec:liftedMeasures}

Let $\lambda$ be the Lebesgue measure on $[0,1]$. Fix $\mu\in\sP$. We call  $\hat\mu\in\Pi(\lambda,\mu)$ a \textit{lift} of $\mu$. The following two ways will be interchangeably used to represent the lifted measure $\hat\mu$.

Let $\hat\theta$ be a measure on $[0,1]\times\R$. Then
\begin{enumerate}
	\item[(i)] $(\hat\theta_{u,\cdot})_{u\in[0,1]}$ denotes the ($\lambda$-a.s. unique) disintegration (w.r.t. the first coordinate) of $\hat\theta$ w.r.t. $\lambda$.
	\item[(ii)] $(\hat\theta_{[0,u],\cdot})_{u\in[0,1]}$ denotes the family of measures defined by
	$$
	\hat\theta_{[0,u],\cdot}(A)=\hat\theta([0,u]\times A)=\int_0^u\hat\theta_{s,\cdot}(A)ds,\quad A\in\sB(\R).
	$$
\end{enumerate}
$\hat\mu_{[0,0],\cdot}$ corresponds to the zero measure on $\R$, while $\hat\mu_{[0,1],\cdot}=\mu$. Furthermore, $\hat\mu_{[0,u],\cdot}\leq\hat\mu_{[0,u'],\cdot}$ for all $u,u'\in[0,1]$ with $u\leq u'$.

We follow Beiglb\"ock and Juillet \cite{BeiglbockJuillet:16s} and call $(\hat\theta_{[0,u],\cdot})_{u\in[0,1]}$ and $(\hat\theta_{u,\cdot})_{u\in[0,1]}$ the  primitive and derivatives curves, respectively. Indeed, $\hat\theta_{u,\cdot}$ can be considered as a derivative of $(\hat\theta_{[0,u],\cdot})_{u\in[0,1]}$ w.r.t. $\sT_1$  (see Beiglb\"{o}ck and Juillet \cite[Section 2.2]{BeiglbockJuillet:16s}). In particular, the set $L\subseteq[0,1]$ of times at which, for any continuous $f$ with linear growth at most, $\int fd\hat\theta_{u,\cdot}=\lim_{h\downarrow 0}(\int fd\hat\theta_{[0,u+h],\cdot}-\int fd\hat\theta_{[0,u],\cdot})/h$ holds, is a (Borel) set of full measure.

Now, in addition to $\hat\mu\in\Pi(\lambda,\mu)$, let $\nu\in\sP$ with $\mu\leq_{cd}\nu$. Then the set of \textit{lifted supermartingale couplings} (or transport plans) is given by
$$
\hat\Pi_S(\hat\mu,\nu):=\left\{\hat\pi\in\Pi(\hat\mu,\nu):\int yd\pi_{(u,x),\cdot}\leq x\textrm{ for $\hat\mu$-a.e. } (u,x)\in[0,1]\times\R\right\},
$$
where $(\hat\pi_{(u,x),\cdot})_{(u,x)\in[0,1]\times\R}$ denotes the disintegration of $\hat\pi$ with respect to $\hat\mu$. Similarly as for a lifted measure $\hat\mu$, we denote the primitive and derivative curves of $\hat\pi$ by $(\hat\pi_{[0,u],\cdot,\cdot})_{u\in[0,1]}$ and $(\hat\pi_{u,\cdot,\cdot})_{u\in[0,1]}$, respectively. Note that both $\hat\pi_{[0,u],\cdot,\cdot}$ and $\hat\pi_{u,\cdot,\cdot}$ are the measures on $\R^2$. Moreover, for any $\hat\pi\in\hat\Pi_S(\hat\mu,\nu)$, the corresponding element of $\Pi_S(\mu,\nu)$ is given by $\pi=\int^1_0\hat\pi_u du$.

For $\hat\mu\in\Pi(\lambda,\mu)$ and $\hat\pi\in\Pi_S(\hat\mu,\nu)$ we have two canonical disintegrations of $\hat\pi$, namely, $\hat\pi(du,dx,dy)=\hat\pi_{u,\cdot,\cdot}(dx,dy)du$ and $\hat\pi(du,dx,dy)=\hat\pi_{(u,x),\cdot}(dy)\hat\mu_{u,\cdot}(dx)du$. In particular, (for $\lambda$-a.e. $u\in[0,1]$ and $\hat\mu_{u,\cdot}$-a.e. $x\in\R$) $\hat\pi_{(u,x),\cdot}$ represents the disintegration of $\hat\pi_{u,\cdot,\cdot}$ with respect to the first marginal $\hat\mu_{u,\cdot}$. Then, if $\hat\nu_{u,\cdot}$ denotes the second marginal of $\hat\pi_{u,\cdot,\cdot}$, we have that $\hat\pi_{u,\cdot,\cdot}\in\Pi_S(\hat\mu_{u,\cdot},\hat\nu_{u,\cdot})$.

{\em Notation:} For $x \in \R$ let $\delta_x$ denote the unit point mass at $x$. For real numbers $c,x,d$ with $c \leq x \leq d$ define the probability measure $\chi_{c,x,d}$ by $\chi_{c,x,d} = \frac{d-x}{d-c}\delta_c + \frac{x-c}{d-c} \delta_d$ if $c<d$ and $\chi_{c,x,d} = \delta_x$ otherwise. (Note that $\chi_{c,x,d}$ has mean $x$ and is the law of a Brownian motion started at $x$ evaluated on the first exit from $(c,d)$.) We extend the definition of $\chi_{c,x,d}$ in the case when one of $\{c,d\}$ takes infinite value. In particular, if $-\infty=c<x\leq d<\infty$ we set $\chi_{c,x,d}=\delta_d$, and similarly, $\chi_{c,x,d}=\delta_c$ whenever $-\infty<c\leq x<d=+\infty$.
\section{The shadow measure and its stability}\label{sec:shadow}
Let $\mu,\nu\in\sM$ with $\mu\leq_{pcd}\nu$. Recall the definition of $\sM^\nu_\mu$, i.e., the set of a target measures of a supermartingale that embeds $\mu$ in $\nu$ (see Section \ref{sec:stochasticOrder} and \eqref{eq:M_nu_mu}). In this section we study the stability properties of the smallest element of $\sM^\nu_\mu$ with respect to $\leq_{cd}$. In what follows, this measure, the so-called \textit{shadow} of $\mu$ in $\nu$, will be denoted by $S^\nu(\mu)$.

Note that, for any $\theta\in\sM^\nu_\mu$ we have that $\overline\theta\leq\overline\mu$. It turns out that we can capture the difference $\overline\mu-\overline{S^\nu(\mu)}$ precisely. For this purpose, for any two measures $\eta,\chi\in\sM$ with $\eta\leq_{pcd}\chi$, we introduce a constant $c_{\eta,\chi}\in[0,\infty)$:
\begin{align}\label{eq:cDef1}
c_{\eta,\chi}:&=\sup_{k\in\R}\{C_\eta(k)-C_\chi(k)\}\\
\label{eq:cDef2} &=\sup_{k\in\R}\{(\chi(\R)-\eta(\R))k-(\overline\chi-\overline\eta)-P_\chi(k)+P_\eta(k)\}.
\end{align}
(The equivalence between \eqref{eq:cDef1} and \eqref{eq:cDef2} is justified by the Put-Call parity.)
\begin{remark}\label{rem:c=0}
Note that, if  $\eta,\chi\in\sM$ with $\eta\leq_{pcd}\chi$, then $c_{\eta,\chi}=0$ if and only if $\eta\leq_{pc}\chi$. Then it follows that $\emptyset\neq\{\theta\in\sM:\mu\leq_{c}\theta\leq\nu\}\subseteq\sM^\nu_\mu$, see Beiglb\"{o}ck at al. \cite[Lemma 4.4]{BHN:20}. In order to verify the claim, note that, by \eqref{eq:cDef1} we have that $c_{\eta,\chi}=0$ if and only if $C_\eta\leq C_\chi$ everywhere (since $\lim_{k\to\infty}C_\eta(k)=\lim_{k\to\infty}C_\chi(k)=0$). On the other hand, Bayraktar et al. \cite[Lemma 3.2]{BDN:21} proved that $C_\eta\leq C_\chi$ (everywhere) if and only if $\eta\leq_{pc}\chi$.
\end{remark}

The next lemma defines, and explicitly determines the shadow measure $S^\nu(\mu)$; see Nutz and Stebegg \cite{NutzStebegg:18} and Bayraktar et al. \cite{BDN:21}.

\begin{lemma}[Shadow measure $S^\nu(\mu)$]\label{lem:shadowDefn}
Let $\mu,\nu\in\sM$ with $\mu\leq_{pcd}\nu$.
\begin{enumerate}
\item There exists a unique measure $S^\nu(\mu)\in\sM^\nu_\mu$ such that $S^\nu(\mu)\leq_{cd}\theta$ for all $\theta\in\sM^\nu_\mu$.
\item The measure $S^\nu(\mu)$ is explicitly determined by its (modified) potential function $P_{S^\nu(\mu)}$, where
\begin{equation}\label{eq:shadowPotential}
P_{S^\nu(\mu)}(k)=P_\nu(k)-(P_\nu-P_\mu)^c(k),\quad k\in\R.
\end{equation}
In particular, $c_{\mu,\nu}=\overline\mu-\overline{S^\nu(\mu)}$.
\item If $\mu=\mu_1+\mu_2$ for some $\mu_1,\mu_2\in\sM$, then
\begin{equation}\label{eq:shadowAssoc}
S^\nu(\mu)=S^\nu(\mu_1)+S^{\nu-S^\nu(\mu_1)}(\mu_2).
\end{equation}
\end{enumerate}
\end{lemma}

We are now ready to present the main result of this section. The following theorem establishes stability of the shadow measure $S^\nu(\mu)$ with respect to initial and target measures $\mu,\nu$.
\begin{theorem}\label{thm:shadowW}
Let $\mu,\mu',\nu,\nu'\in\sM$ with $\mu(\R)=\mu'(\R)$ and $\nu(\R)=\nu'(\R)$. Suppose $\mu\leq_{pcd}\nu$ and $\mu'\leq_{pcd}\nu'$. The following relation holds
\begin{equation}\label{eq:shadowW}
W(S^\nu(\mu),S^{\nu'}(\mu'))\leq W(\mu,\mu')+2W(\nu,\nu').
\end{equation}
\end{theorem}

The proof of Theorem \ref{thm:shadowW} relies on two auxiliary propositions and the following \textit{up} and \textit{down} measures.

For $\eta,\eta'\in\sM$ with $\eta(\R)=\eta'(\R)$ let Up$(\eta,\eta')$, Down$(\eta,\eta')\in\sM$ be given by
\begin{align}
\textrm{Up}(\eta,\eta')(B)&=\lambda\lvert_{[0,\eta(\R)=\eta'(\R)]}(\{u\in\R:(G_\eta(u)\vee G_{\eta'}(u))\in B\}),&\quad B\in\sB(\R),\label{eq:Up}\\
\textrm{Down}(\eta,\eta')(B)&=\lambda\lvert_{[0,\eta(\R)=\eta'(\R)]}(\{u\in\R:(G_\eta(u)\wedge G_{\eta'}(u))\in B\}),&\quad B\in\sB(\R).\label{eq:Down}
\end{align}
Observe that  Up$(\eta,\eta')(\R)=$ Down$(\eta,\eta')(\R)=\eta(\R)=\eta'(\R)$.

\begin{proposition}\label{prop:shadowW1}
Let $\mu,\mu',\nu\in\sM$ be such that $\mu(\R)=\mu'(\R)$, $\mu\leq_{pcd}\nu$ and $\mu'\leq_{pcd}\nu$. We have
	$$
	W(S^\nu(\mu),S^\nu(\mu'))\leq W(\mu,\mu').
	$$
\end{proposition}
\begin{proposition}\label{prop:shadowW2}
Let $\mu,\nu,\nu'\in\sM$ be such that $\mu\leq_{pcd}\nu$, $\mu\leq_{pcd}\nu'$ and $\nu(\R)=\nu'(\R)$. Then
	$$
	W(S^\nu(\mu),S^{\nu'}(\mu))\leq 2 W(\nu,\nu').
	$$
\end{proposition}
\begin{proof}[Proof of Theorem \ref{thm:shadowW}]
First, by Lemma \ref{lem:UpDownSto}, we have that  Down$(\nu,\nu')\leq_{sto}\chi$ for $\chi\in\{\nu,\nu'\}$, while by Lemma \ref{lem:UpDownPCD}, $\eta\leq_{pcd}$ Down$(\nu,\nu')$ for $\eta\in\{\mu,\mu'\}$. Then
\begin{align*}
&W(S^\nu(\mu),S^{\nu'}(\mu'))\nonumber\\
&\leq W(S^\nu(\mu),S^{\textrm{Down}(\nu,\nu')}(\mu)) + W(S^{\textrm{Down}(\nu,\nu')}(\mu),S^{\textrm{Down}(\nu,\nu')}(\mu'))+W(S^{\textrm{Down}(\nu,\nu')}(\mu'),S^{\nu'}(\mu'))\\
&\leq 2W(\nu,\textrm{Down}(\nu,\nu')) + W(\mu,\mu')+2W(\nu',\textrm{Down}(\nu,\nu')),
\end{align*}
where the first inequality follows from the triangle inequality (applied twice), while for the second one we used Propositions \ref{prop:shadowW1} and  \ref{prop:shadowW2}. We finish the proof by observing that $W(\nu,\textrm{Down}(\nu,\nu')) +W(\nu',\textrm{Down}(\nu,\nu'))=W(\nu,\nu')$, see Lemma \ref{lem:UpDownW}. 
\end{proof}
\begin{remark}\label{rem:wassertstein}
Theorem \ref{thm:shadowW} is a supermartingale generalisation of Juillet \cite[Theorem 2.31]{Juillet:16}. While our proof uses similar structure, it is simpler and more direct.

For example, in several instances the proof relies on the explicit representation of the shadow measure $S^\nu(\mu)$. In the martingale case, Juillet \cite{Juillet:16} first uses the representation of $S^\nu(\mu)$ when $\mu$ is atomic, and then obtains general statements by approximation. In our case, and directly for general $\mu$, we use the representation of $S^\nu(\mu)$ via (modified) potential function $P_{S^\nu(\mu)}$ (see Lemma \ref{lem:shadowDefn}), and thus bypass the approximation step.

Furthermore, a supermartingale has a natural direction (i.e., a tendency to decrease), which translates to convenient relations between measures in terms of stochastic orders. To illustrate this, let $\mu,\nu,\nu'$ be as in the statement of Theorem\ref{thm:shadowW}, i.e., $\mu\leq_{pcd}\nu$ and $\mu\leq_{pcd}\nu'$. Then $\mu\leq_{pcd}\textrm{\normalfont{Down}}(\nu,\nu')$, see Lemma \ref{lem:UpDownPCD}. On the other hand, if $\mu\leq_{pc}\nu$ and $\mu\leq_{pc}\nu'$ then, in general, $\mu\leq_{pc}\textrm{\normalfont{Down}}(\nu,\nu')$ does not hold, and thus in the martingale case Juillet \cite{Juillet:16} needs additional arguments.
\end{remark}
\subsection{Proofs of Propositions \ref{prop:shadowW1} and \ref{prop:shadowW2}}
The proofs of both propositions rely on the following important lemma.
\begin{lemma}\label{lem:1stOrder}
Let $\mu,\mu',\nu\in\sM$ with $\mu(\R)=\mu'(\R)$, $\mu\leq_{pcd}\nu$ and $\mu'\leq_{pcd}\nu$. If $\mu\leq_{sto}\mu'$, then $S^\nu(\mu)\leq_{sto}S^\nu(\mu')$.
\end{lemma}
\begin{proof}
	Note that $S^\nu(\mu)\leq_{sto}S^\nu(\mu')$ is equivalent to $F_{S^\nu(\mu)}\geq F_{S^\nu(\mu')}$ on $\R$. On the other hand, for any $\eta\in\sM$, $F_\eta(x)=\eta((-\infty,x])=(P_\eta)'_+(x)$ for each $x\in\R$. Hence it is enough to show that
	\begin{equation}\label{eq:potentialDerivative}
	(P_{S^\nu(\mu)})'_+(x)\geq (P_{S^\nu(\mu')})'_+(x),\quad x\in\R.
	\end{equation}
	By Lemma \ref{lem:shadowDefn} (see \eqref{eq:shadowPotential}), \eqref{eq:potentialDerivative} will be established if we can show that
	$$
	[(P_\nu-P_\mu)^c]'_+(x)\leq[(P_\nu-P_{\mu'})^c]'_+(x),\quad x\in\R.
	$$
	But since $\mu\leq_{sto}\mu'$,
	$$
	(P_\nu-P_\mu)'_+(x)\leq(P_\nu-P_{\mu'})'_+(x),\quad x\in \R.
	$$
	Then an application of Lemma \ref{lem:hullDerivative} with $f=(P_\nu-P_\mu)$ and $g=(P_\nu-P_{\mu'})$ completes the proof.
\end{proof}	
\begin{proof}[Proof of Proposition \ref{prop:shadowW1}]
To ease the notation, set $\tilde\mu=$Up$(\mu,\mu')$. From Lemma \ref{lem:UpDownSto} we have that $\mu\leq_{sto}\tilde\mu$  and $\mu'\leq_{sto}\tilde\mu$, while Lemma \ref{lem:Up1marginals} ensures that $\tilde\mu\leq_{pcd}\nu$. This permits us to apply Lemma \ref{lem:1stOrder}, from which we conclude that $S^\nu(\mu)\leq_{sto}S^\nu(\tilde\mu)$ and $S^\nu(\mu')\leq_{sto}S^\nu(\tilde\mu)$.
	 
Using the above observations we have that
	\begin{align}
	W(\mu,\mu')&=W(\mu,\tilde\mu)+W(\mu',\tilde\mu)\nonumber\\
	&=\overline{\tilde\mu}-\overline\mu + \overline{\tilde\mu}-\overline{\mu'}\nonumber\\
	&=W(S^\nu(\mu),S^\nu(\tilde\mu))+W(S^\nu(\mu'),S^\nu(\tilde\mu))\label{eq:prop21proof}\\
	&+\{\overline{\tilde\mu}-\overline{S^\nu(\tilde\mu)}-(\overline\mu-\overline{S^\nu(\mu)})\}\nonumber\\
	&+\{\overline{\tilde\mu}-\overline{S^\nu(\tilde\mu)}-(\overline{\mu'}-\overline{S^\nu(\mu')})\}\nonumber,
	\end{align}
	where the first equality follows from Lemma \ref{lem:UpDownW}, while the other two use the fact (see Section \ref{sec:measureSpaces}) that $W(\eta,\chi)=\overline\chi-\overline\eta$ whenever $\eta,\chi\in\sM$ are such that $\eta(\R)=\chi(\R)$ and $\eta\leq_{sto}\chi$.

We claim that the last two summands in \eqref{eq:prop21proof} are non-negative. We will only verify this for $\mu$, i.e., we consider the penultimate summand (the case for $\mu'$ is identical). By Lemma \ref{lem:shadowDefn} we have that
	$$
	\overline{\tilde\mu}-\overline{S^\nu(\tilde\mu)}-(\overline\mu-\overline{S^\nu(\mu)})=c_{\tilde\mu,\nu}-c_{\mu,\nu}.
	$$
	But since, for each $k\in\R$, $s\mapsto (s-k)^+$ is non-decreasing, and also $\mu\leq_{sto}\tilde\mu$, we have that $C_{\tilde\mu}(k)\geq C_{\mu}(k)$, and it follows that
	$$
	c_{\tilde\mu,\nu}-c_{\mu,\nu}=\sup_{k\in\R}\{C_{\tilde\mu}(k)-C_\nu(k)\}-\sup_{k\in\R}\{C_{\mu}(k)-C_\nu(k)\}\geq0.
	$$

Finally,
$$
	W(\mu,\mu')\geq W(S^\nu(\mu),S^\nu(\tilde\mu))+W(S^\nu(\mu'),S^\nu(\tilde\mu))\geq W(S^\nu(\mu),S^\nu(\mu')),
	$$
	where the last inequality follows from the triangle inequality.
\end{proof}

The proof of Proposition \ref{prop:shadowW2} will need one additional result.
\begin{lemma}\label{lem:shadowConvergence}
Let $\mu,\nu\in\sM$ with $\mu\leq_{pcd}\nu$, and let $(\mu_n)_{n\geq1}$ be a sequence of measures in $\sM$ increasing in convex order and such that $\mu_n\leq_c\mu$.

Then $\mu_n\to\mu_\infty$ and $S^\nu(\mu_n)\to S_\infty$ for some $\mu_\infty, S_\infty\in\sM$. In particular, $S_\infty=S^\nu(\mu_\infty)$.
\end{lemma}
\begin{proof}
It is well-known that $\mu_n\xrightarrow{w}\mu_\infty$ (increasingly with respect to convex order) for some $\mu_\infty\in\sM$ if and only if $U_{\mu_n}\downarrow U_{\mu_\infty}$ pointwise (see, for example, Chacon \cite{Chacon:77}). In this case, the first moments also converge (i.e., $U_{\mu_n}(0)\downarrow U_{\mu_\infty}(0)$), and therefore $\mu_n\xrightarrow{w}\mu_\infty$ is equivalent to $\mu_n\to\mu$.

Since $\mu_n\leq_c\mu_{n+1}\leq_c\mu$, $U_{\mu_n}(k)\geq U_{\mu_{n+1}}(k)\geq U_\mu(k)$ and hence $\lim_{n\to\infty}U_{\mu_n}(k)$ exists for each $k\in\R$. It is easy to see that $U_\infty:=\lim_{n\to\infty}U_{\mu_n}$ is concave and with the same asymptotic behaviour as $U_{\mu_n}$ and $U_\mu$. It follows that $U_\infty=U_{\mu_\infty}$ for some $\mu_\infty\in\sM$ with the same mean and mass as $\mu$. We conclude that $\mu_n\to\mu_\infty$. Furthermore, since $U_{\mu_\infty}=U_{\infty}\geq U_\mu$, we have that $\mu_\infty\leq_c\mu$, and therefore $\mu_\infty\leq_{pcd}\nu$, so that the shadow measure $S^\nu(\mu_\infty)$ is well-defined.

Now note that, for $n\geq1$,
\begin{align*}
P_{S^\nu(\mu_n)}=P_\nu-(P_\nu-P_{\mu_n})^c&=P_{S^\nu(\mu_{n+1})}+(P_\nu-P_{\mu_{n+1}})^c-(P_\nu-P_{\mu_n})^c\\
&=P_{S^\nu(\mu_\infty)}+(P_\nu-P_{\mu_\infty})^c-(P_\nu-P_{\mu_{n}})^c.
\end{align*}
Therefore, since $P_{\mu_n}\leq P_{\mu_{n+1}}\leq P_{\mu_\infty}$, we have that $P_{S^\nu(\mu_n)}\leq P_{S^\nu(\mu_{n+1})}\leq P_{S^\nu(\mu_\infty)}$. It follows that, for all $n\geq 1$, $S^\nu(\mu_n)\leq_{cd}S^\nu(\mu_{n+1})\leq_{cd} S^\nu(\mu_\infty)$ and therefore $\overline{S^\nu(\mu_n)}\geq\overline{S^\nu(\mu_{n+1})}\geq\overline{S^\nu(\mu_\infty)}$. In particular, $\lim_{n\to\infty}\overline{S^\nu(\mu_n)}$ exists and (since $\nu$ is integrable) it is finite. Furthermore, applying Lemma \ref{lem:converge_hull} with $f_n=(P_\nu-P_{\mu_n})$ and $f=(P_\nu-P_{\mu_\infty})$ we have that $P_{S^\nu(\mu_n)}\uparrow P_{S^\nu(\mu_\infty)}$ pointwise, as $n\to\infty$.

Recall that $U_\eta=-C_\eta-P_\eta$ for $\eta\in\sM$. Using the Put-Call parity we further have that $U_\eta(k)=\eta(\R)k-\overline\eta-2P_\eta(k)$ for each $k\in\R$. Then, since $S^\nu(\mu_n)(\R)=\mu_n(\R)=\mu_\infty(\R)=S^\nu(\mu_\infty)(\R)$, we have that, for each $k\in\R$,
\begin{align*}
\lim_{n\to\infty}U_{S^\nu(\mu_n)}(k)&=S^\nu(\mu_\infty)(\R)k-\lim_{n\to\infty}\overline{S^\nu(\mu_n)}-2P_{S^\nu(\mu_\infty)}(k)\\&=U_{S^\nu(\mu_\infty)}(k)-\left\{\lim_{n\to\infty}\overline{S^\nu(\mu_n)}-\overline{S^\nu(\mu_\infty)}\right\}.
\end{align*}
It follows that $\lim_{n\to\infty}U_{S^\nu(\mu_n)}(k)$ exists for all $k\in\R$, and by Chacon \cite[Lemma 2.6]{Chacon:77}, $S^\nu(\mu_n)\xrightarrow{w} S_\infty$ for some $S_\infty\in\sM$.

Finally, let $U_\infty:=\lim_{n\to\infty}U_{S^\nu(\mu_n)}$ and $C:=\lim_{n\to\infty}\overline{S^\nu(\mu_n)}-\overline{S^\nu(\mu_\infty)} \geq0$. Then applying Chacon \cite[Lemma 2.5]{Chacon:77} we have that
$$
U_{S_\infty}=U_{\infty} +C=U_{S^\nu(\mu_\infty)}.
$$
Since the potential functions uniquely identify the underlying measures, it follows that $$S^\nu(\mu_n)\xrightarrow{w}S_\infty=S^\nu(\mu_\infty).$$

It is left to show that $C=0$ (from which we can conclude that the first moments of $S^\nu(\mu_n)$ converge to the first moment of $S^\nu(\mu_\infty)$) and therefore $S^\nu(\mu_n)\xrightarrow{w} S_\infty$ is equivalent to $S^\nu(\mu_n)\to S_\infty$. By Billingsley \cite[Theorem 3.5]{Billingsley:2013} it is enough to show that $\left\{S^\nu(\mu_n)\right\}_{n\geq 1}$ is uniformly integrable. Note that $S^\nu(\mu_n)(\R)=\mu(\R)$ for all $n\geq1$. Let $\underline{\theta},\overline{\theta}\in\sM$ be such that $\underline\theta(\R)=\overline\theta(\R)=\mu(\R)$, $\underline\theta\leq\nu$, $\overline\theta\leq\nu$ and $\underline\theta\leq_{sto}\eta\leq_{sto}\overline\theta$ for all $\eta\in\sM$ with $\eta(\R)=\mu(\R)$ and $\eta\leq\nu$. Then $\underline\theta,\overline\theta$ are the `left-most' and `right-most' measures (of mass $\mu(\R)$) within $\nu$ (i.e., $\underline\theta$ is the restriction of $\nu$ between 0-th and $\mu(\R)$-th quantiles (see \eqref{eq:left-most}), while $\overline\theta$ is the restriction of $\nu$ between $(\nu(\R)-\mu(\R))$-th and $\nu(\R)$-th quantiles), respectively, and we have that 
\begin{align*}
0&\leq\lim_{K\to\infty}\sup_{n\geq 1}\int_{(-\infty,-K]\cup[K,+\infty)}\lvert x\lvert dS^\nu(\mu_n)(x)\\
&\leq \left(\lim_{K\to+\infty}\int_{(-\infty,-K]}\lvert x\lvert d\underline\theta(x)\right)+\left(\lim_{K\to+\infty}\int_{[K,+\infty)}\lvert x\lvert d\overline\theta(x)\right)\\
=&\lim_{K\to\infty}\int_{(-\infty,-K]\cup[K,+\infty)}\lvert x\lvert d\nu(x)= 0,
\end{align*}
where for the equalities we use that $\underline\theta$ and $\overline\theta$ are the restrictions of $\nu$, and that $\nu$ is integrable, respectively.
\end{proof}

\begin{proof}[Proof of Proposition \ref{prop:shadowW2}] Let $\tilde\nu=$Down$(\nu,\nu')$. Note that, by Lemma \ref{lem:UpDownSto}, $\tilde\nu\leq_{sto}\nu$ and $\tilde\nu\leq_{sto}\nu'$, while Lemma \ref{lem:UpDownPCD} ensures that $\mu\leq_{pcd}\tilde\nu$. From the triangle inequality we deduce that $W(S^\nu(\mu),S^{\nu'}(\mu))\leq	W(S^\nu(\mu),S^{\tilde\nu}(\mu))+	W(S^{\tilde\nu}(\mu),S^{\nu'}(\mu))$.	
	
	Suppose the claim of Proposition \ref{prop:shadowW2} is true for $\nu'\leq_{sto}\nu$. Then $W(S^\nu(\mu),S^{\tilde\nu}(\mu))\leq 2 W(\nu,\tilde\nu)$ and $W(S^{\tilde\nu}(\mu),S^{\nu'}(\mu))\leq 2 W(\nu',\tilde\nu)$, and then using Lemma \ref{lem:UpDownW} we obtain $$W(S^\nu(\mu),S^{\nu'}(\mu))\leq 2 W(\nu,\nu'),$$	as required.
	
	It is left to show that the claim of Proposition \ref{prop:shadowW2} holds under an assumption that $\nu'\leq_{sto}\nu$. Because of Lemma \ref{lem:shadowConvergence} we can further assume that $\mu$ is of type $\sum^n_{i=1}\alpha_i\delta_{x_i}$ (the general $\mu$ is then approximated by a sequence of atomic measures increasing in convex order).

Then, by Lemma \ref{lem:constructionMU}, there exists $\mu'\in\sM$ with $\mu'(\R)=\mu(\R)$, and such that $\mu\leq_{sto}\mu'$, $S^{\nu'}(\mu)\leq_{sto}S^\nu(\mu')$ and $W(S^{\nu'}(\mu),S^{\nu}(\mu'))\leq W(\nu,{\nu'})$ . From Lemma \ref{lem:1stOrder} we also have that $S^\nu(\mu)\leq_{sto} S^\nu(\mu')$, and it follows that
	\begin{align}
	W(S^\nu(\mu),S^{\nu'}(\mu))&\leq W(S^{\nu'}(\mu),S^{\nu}(\mu'))+W(S^{\nu}(\mu'),S^{\nu}(\mu))\nonumber\\
	&=\overline{S^{\nu}(\mu')}-\overline{S^{\nu'}(\mu)} + \overline{S^{\nu}(\mu')}-\overline{S^{\nu}(\mu)}\label{eq:prop2.2}\\
&\leq 2 W(\nu',\nu)+\overline{S^{\nu'}(\mu)}-\overline{S^{\nu}(\mu)},\nonumber
\end{align}
where we used the triangle inequality and the fact (see Section \ref{sec:Wasserstein} and \eqref{eq:WassQuantile}) that $W(\eta,\chi)=\int^{\eta(\R)}_0G_\chi(u)-G_\eta(u)du=\overline\chi-\overline\eta$ whenever $\eta(\R)=\chi(\R)$ and $\eta\leq_{sto}\chi$. 

Hence the claim of the proposition holds if we can show that $\overline{S^{\nu'}(\mu)}-\overline{S^{\nu}(\mu)}\leq 0$. However, by Lemma \ref{lem:shadowDefn},
	\begin{align*}
\overline{S^{\nu'}(\mu)}-\overline{S^{\nu}(\mu)}&=\overline\mu-c_{\mu,\nu'}+c_{\mu,\nu}-\overline\mu\\
&=\sup_{k\in\R}\{C_\mu(k)-C_\nu(k)\}-\sup_{k\in\R}\{C_\mu(k)-C_{\nu'}(k)\}\leq 0,
	\end{align*}
where the inequality follows from the fact that $C_{\nu'}\leq C_{\nu}$, since $\nu'\leq_{sto}\nu$ and, for each $k\in\R$, $s\mapsto(s-k)^+$ is non-decreasing. 
\end{proof}

In the proof of Proposition \ref{prop:shadowW2} we used the following lemma, which can be proved using a construction provided in the first part of the proof of Juillet \cite[Proposition 2.36]{Juillet:16}. We include the proof for the convenience of the reader.
\begin{lemma}\label{lem:constructionMU}
Let $\mu,\nu,\nu'\in\sM$ be such that $\nu(\R)=\nu'(\R)$, $\nu'\leq_{sto}\nu$, $\mu\leq_{pcd}\nu$ and $\mu\leq_{pcd}\nu'$. Suppose $\mu$ is of the form $\sum^n_{i=1}\alpha_i\delta_{x_i}$ for some $n\in\mathbb{N}$. Then there exists $\mu'\in\sM$ such that $\mu'(\R)=\mu(\R)$, $\mu\leq_{sto}\mu'$, $S^{\nu'}(\mu)\leq_{sto}S^\nu(\mu')$ and $W(S^{\nu'}(\mu),S^{\nu}(\mu'))\leq W(\nu,{\nu'})$.
\end{lemma}
\begin{proof}
The construction of $\mu'$ relies on the following fact: if $\mu$ is as in the statement and $n=1$, then $S^{\nu'}(\mu)$ is a restriction of $\nu'$ between two quantiles (see Nutz and Stebegg \cite[Lemma 6.3]{NutzStebegg:18}). More generally (when $n>1$), there exists a sequence of sets $J_1\subseteq\dots\subseteq J_n\subseteq(0,\nu(\R)=\nu'(\R)]$ such that, for any $k\in\{1,\dots,n\}$, $S^{\nu'}(\sum^k_{i=1}\alpha_i\delta_{x_i})(B)=\lambda\lvert_{[0,\nu'(\R)]}(\{u\in J_k:G_{\nu'}(u)\in B\})$, for all Borel subsets $B$ of $\R$.
	
	Now introduce $\mu'=\sum^n_{i=1}\alpha_i\delta_{x_i'}$, where $x'_i$ is the barycenter of a measure $\eta_i$, defined by $\eta_i(B)=\lambda\lvert_{[0,\nu(\R)]}(\{u\in J_i\setminus J_{i-1}:G_{\nu}(u)\in B\})$. Since $\nu'\leq_{sto}\nu$, $G_{\nu'}\leq G_\nu$ and therefore $x_i\leq x'_i$. It follows that $\mu\leq_{sto}\mu'$ and (since $\mu\leq_{pcd}\nu$) $\mu'\leq_{pcd}\nu$. It is then easy to see that the shadow of $\mu'$ in $\nu$ is given by $S^{\nu}(\mu')(B)=\lambda\lvert_{[0,\nu(\R)]}(\{u\in J_n:G_{\nu}(u)\in B\})$ and, in particular, $S^{\nu'}(\mu)\leq_{sto} S^\nu(\mu')$. 

Finally,
\begin{align*}
W(S^{\nu'}(\mu),S^{\nu}(\mu'))=\overline{S^{\nu}(\mu')}-\overline{S^{\nu'}(\mu)}&=\int_{J_n}[G_\nu(u)-G_{\nu'}(u)]du\\
&\leq \int_0^{\nu(\R)=\nu'(\R)}[G_\nu(u)-G_{\nu'}(u)]du= W(\nu,{\nu'}).
\end{align*}
\end{proof}

\section{Existence of lifted shadow couplings}\label{sec:mainExistence}
In this section we show that the shadow measure allows to construct a large family of (lifted) supermartingale couplings.

Let $\sF(\R)$ be the space of closed subsets of $\R$. The space $\sF(\R)$ is endowed with the coarsest topology such that $F\in\sF(\R)\mapsto d(x,F)$ is continuous for every $x\in\R$ (see Kellerer \cite[Section 2.1]{kellerer73}). We write $T\in\sI$ if $T\in\sF(\R)$ and $\inf \{k\in T\}=-\infty$.

\begin{defn}
	Let $T\in\sI$. For every $x\in\R$, let $x^-_T:=\sup\{k\in T\cap(-\infty, x]\}$ and $x^+_T:=\inf\{k\in T\cap[x,\infty)\}$ with $\inf\emptyset=\infty$. The Kellerer dilation is given by
	\begin{align*}
	D_T(x,\cdot)=\begin{cases}
	\delta_x& \textrm{if }x\in T;\\
	\frac{x^+_T-x}{x^+_T-x^-_T}\delta_{x^-_T}+\frac{x-x^-_T}{x^+_T-x^-_T}\delta_{x^+_T}&\textrm{otherwise}.
	\end{cases}
	\end{align*}
	(If $x\notin T$ and $x^+_T=\infty$, then $D_T(x,\cdot)=\delta_{x^-_T}$.)
\end{defn}
If $\mu\in\sP$, then the \textit{hitting projection} of $\mu$ in $T$ is a measure $\mu D_T(B)=\int_\R D_T(x,B)d\mu(x)$, and the \textit{hitting coupling} of $\mu$ and $\mu D_T$ is defined by $\pi_{\mu,T}(A\times B)=\int_A D_T(x,B)d\mu(x)$.

Note that if $T$ is not an element of $\sI$, but $\inf \{k\in T\}\leq \inf\{x\in supp(\mu)\}$, then the kernel $D_T$ still makes sense $\mu$-a.s.

\begin{remark}\label{rem:suermartingaleDilation}
The kernel $D_T$ is slightly different from the original dilation introduced in Kellerer \cite[Definition 16]{kellerer73}. The difference lies in the definition of $\sI$. Kellerer \cite{kellerer73}  considers martingales only, and therefore writes $T\in\sI$ if $T\in\sF(\R)$, $\inf \{k\in T\}=-\infty$ and $\sup \{k\in T\}=\infty$. The last condition, however, cannot be guaranteed in the supermartingale setting. (For example, consider a supermartingale with the starting law $\mu$ supported on $(0,\infty)$ and the target law $\nu$ supported on $(-\infty,0)$, and take $x\in 
supp(\mu)$ together with $T=supp(\nu)$.) Therefore in the supermartingale case we need to explicitly deal with a situation when $\{k\in T\cap[x,\infty)\}$ is empty, $x\in\R$.
\end{remark}
Before stating the main result of this section, we present a useful result regarding the dilation $D_T$. The proof is postponed until Appendix \ref{sec:AppProofs}.
\begin{lemma}\label{lem:continuousProjection}
	Fix $T\in\sF(\R)$ and $\mu\in\sP$ with $\inf \{k\in T\}\leq\inf\{k\in supp(\mu)\}$.
	\begin{enumerate}
		\item The hitting coupling $\pi_{\mu,T}$ is the unique element of $\Pi_S(\mu,\mu D_T)$.
		\item Let $(\mu_n)_{n\geq 1}$ be a sequence in $\sP$ with $\mu_n(\R)=\mu(\R)$ for all $n\geq1$. If $\mu_n\to\mu$ then $\mu_nD_T\to\mu D_T$.
	\end{enumerate}
\end{lemma}	
With the help of dilation $D_T$ we can now formulate the main result of this section. Recall that $\lambda$ denotes the Lebesgue measure on $[0,1]$; see also Section \ref{sec:liftedMeasures} for the definitions of primitive and derivative curves of a lifted measure.
\begin{theorem}\label{thm:shadowCouplings}
Let $\mu,\nu\in\sP$ with $\mu\leq_{cd}\nu$, and let $\hat\mu\in\Pi(\lambda,\mu)$. Then there exists a unique element $\hat\pi\in\Pi_S(\hat\mu,\nu)$, the lifted shadow coupling of $\hat\mu$ and $\nu$, such that for every $u\in[0,1]$, the first and second marginals of $\hat\pi_{[0,u],\cdot.\cdot}$ are $\hat\mu_{[0,u],\cdot}$ and $\hat\nu_{[0,u],\cdot}:=S^\nu(\hat\mu_{[0,u],\cdot})$, respectively. If we denote by $\hat\nu_{u,\cdot}$ the derivative of $(\hat\nu_{[0,u],\cdot})_{u\in[0,1]}$ at $u$ (whenever it exists), we have moreover $\hat\nu_{u,\cdot}=\hat\mu_{u,\cdot} D_{T(u)}$, where $T(u):=supp(\nu-\hat\nu_{[0,u],\cdot})$.
\end{theorem}
The proof of Theorem \ref{rem:suermartingaleDilation} relies on the following result.
\begin{proposition}\label{prop:Dilation}
Let $\mu,\nu\in\sP$ with $\mu\leq_{cd}\nu$, and let $\hat\mu\in\Pi(\lambda,\mu)$. Suppose $u_0\in[0,1)$ is such that $u\mapsto\hat\mu_{[0,u],\cdot}$ has a right derivative at $u_0$, and let $\hat\nu_{[0,u],\cdot}=S^\nu(\hat\mu_{[0,u],\cdot})$. Then a right derivative of $(\hat\nu_{[0,u],\cdot})_{u\in[0,1]}$ at $u_0$ exists and is given by $\hat\mu_{u_0,\cdot}D_T$, where $T=supp(\nu-\hat\nu_{[0,u_0],\cdot})$. Furthermore, $\inf \{k\in T\}\leq \inf\{k\in supp(\hat\mu_{u_0,\cdot})\}$.
\end{proposition}

\begin{proof}[Proof of Theorem \ref{thm:shadowCouplings}]
	Using Proposition \ref{prop:Dilation}, the proof can be obtained by the arguments of Beiglb\"{o}ck and Juillet \cite[Theorem 2.9]{BeiglbockJuillet:16}. Nevertheless, we sketch the proof for completeness.

Let $\hat\mu_{[0,u],\cdot}$, $\hat\nu_{[0,u],\cdot}$ and $T(u)$ be as in the statement. Then using  Proposition \ref{prop:Dilation}, and by setting $\hat\nu_{u,\cdot}=\hat\mu_{u,\cdot}D_{T(u)}$, we can define $\hat\pi_{u,\cdot,\cdot}=\pi_{\hat\mu_{u,\cdot},T(u)}\in\Pi_S(\hat\mu_{u,\cdot},\hat\nu_{u,\cdot})$ for almost every $u\in[0,1]$ (recall that $\pi_{\hat\mu_{u,\cdot},T(u)}$ denotes the hitting coupling), and then the associated $\hat\pi\in\hat\Pi_S(\hat\mu,\nu)$ and $\hat\pi_{[0,u],\cdot,\cdot}=\int^u_0\hat\pi_{t,\cdot,\cdot}dt$.

On the other hand, if $\hat\pi\in\hat\Pi_S(\hat\mu,\nu)$ is such that the marginals of $\hat\pi_{[0,u],\cdot,\cdot}$ are $\hat\mu_{[0,u],\cdot}$ and $\hat\nu_{[0,u],\cdot}=S^\nu(\hat\mu_{[0,u],\cdot})$, then at points $u$ where the derivatives of $\hat\mu_{[0,u],\cdot}$, $\hat\nu_{[0,u],\cdot}$ and $\hat\pi_{[0,u],\cdot,\cdot}$ exist we have that $\hat\pi_{u,\cdot,\cdot}\in\Pi_S(\hat\mu_{u,\cdot},\hat\nu_{u,\cdot})$. But by Proposition \ref{prop:Dilation} we necessarily have that $\hat\nu_{u,\cdot}=\hat\mu_{u,\cdot}D_{T(u)}$. Then the uniqueness part of Lemma \ref{lem:continuousProjection} completes the proof.
\end{proof}
It is left to prove Proposition \ref{prop:Dilation}. We will need the following lemma.

\begin{lemma}\label{lem:convergence}
	Let $(H_n)_{n\geq1}$ be a sequence of positive numbers tending to infinity, $(\eta_n)_{n\geq1}$ a sequence in $\sP$ converging to $\eta\in\sP$, and $v\in\sM$. Assume $\eta_n\leq_{pcd}H_nv$ for every $n\geq1$. Then, setting $T=supp(v)$, it holds $\inf\{k\in T\}\leq\inf\{k\in supp(\eta)\}$ and $S^{H_nv}(\eta_n)\to\eta D_T$ in $\sP$.
\end{lemma}
\begin{proof}
	The proof follows the same arguments as in the martingale case of Beiglb\"{o}ck and Juillet \cite[Lemma 2.8]{BeiglbockJuillet:16s}. The main difference lies in the assumption $\eta_n\leq_{pcd} H_n v$ (and not $\leq_{pc}$) and the definition of $D_T$. Therefore we will only highlight the necessary modifications.
	
	First, since $\eta_n\leq_{pcd}H_n v$ we have that $\inf \{k\in T\}\leq \inf\{k\in supp(\eta_n)\}$ for every $n\geq1$, and therefore $\eta_n([\inf T,\infty))=1$. Letting $n\to\infty$ we find $supp(\eta)\subseteq[\inf \{k\in T\},\infty]$. 

		1. First suppose that $\eta_n=\delta_x$ for all $n\geq1$. The proof in the martingale case relies on the fact that (when $\delta_x\leq_{pc} H_n v$) the martingale shadow $S^{H_nv}(\delta_x)$ is supported on an interval. However, the same is true in the supermartingale case, see Nutz and Stebegg \cite[Lemma 6.3]{NutzStebegg:18}.
		
		2. Now suppose that $\eta_n=\eta=\sum^n_{k=1}m_k\delta_{x_k}$ for all $n\geq 1$. The proof in the martingale case relies on Step 1, associativity of the (martingale) shadow measure and the induction argument. Since the associativity also holds in the supermartingale case (recall Lemma \ref{lem:shadowDefn}), the result follows.

3. In this step the result is established for a constant sequence $\eta_n=\eta$. To achieve this, approximate a general measure $\eta$ by atomic measures $(\eta_k)_{k\geq 1}$ with $\eta_k\leq _c\eta$ and such that $\eta_k\to\eta$ in $\sP$ as $k\to\infty$. Note that all the measures $(\eta_k)_{k\geq 1},\eta$ have the same total mass and mean, which will enable us to use Lemma \ref{lem:continuousProjection}.

To establish the claim (and similarly as in the martingale case) we use Step 2 and the following two facts: first, by Theorem \ref{thm:shadowW} we have that
\begin{equation}\label{Wdistance}
W(S^{H_n v}(\eta_k),S^{H_n v}(\eta))\leq W(\eta_k,\eta)
\end{equation}
(which converges to zero uniformly in $n$ as $k$ goes to infinity), and second, $\eta_k D_T\to\eta D_T$ which is guaranteed by Lemma \ref{lem:continuousProjection}.

4. If $\eta_n$ is a non-constant sequence, then note that
\begin{align*}
W(S^{H_n v}(\eta_n),\eta D_T)&\leq W(S^{H_n v}(\eta_n),S^{H_n v}(\eta)) + W(S^{H_nv}(\eta),\eta D_T)\\
&\leq W(\eta_n,\eta)+ W(S^{H_nv}(\eta),\eta D_T)
\end{align*}
where the second inequality follows from Theorem \ref{thm:shadowW}. Since $\eta_n\to\eta$, $W(\eta_n,\eta)\to0$, while Step 3 ensures that $W(S^{H_nv}(\eta),\eta D_T)\to0$.
\end{proof} 

Finally we can prove Proposition \ref{prop:Dilation}; the arguments are identical to those of Beiglb\"ock and Juillet \cite[Proposition 2.7]{BeiglbockJuillet:16s} and thus we only give a sketch.
\begin{proof}[Proof of Proposition \ref{prop:Dilation}]
	For $h>0$ consider $h^{-1}(\hat\nu_{[0,u_0+h],\cdot}-\hat\nu_{[0,u_0],\cdot})=h^{-1}(S^\nu(\hat\mu_{[0,u_0+h],\cdot})-S^\nu(\hat\mu_{[0,u_0],\cdot}))=:\sigma_h\in\sP$. By associativity of the shadow measure and by using an appropriate scaling of measures we have that $\sigma_h=S^{h^{-1}(\nu-\hat\nu_{[0,u_0],\cdot})}(h^{-1}(\hat\mu_{[0,u_0+h],\cdot}-\hat\mu_{[0,u_0],\cdot}))$. Replacing $h$ by a sequence $(h_n)_{n\geq1}$ (of positive numbers decreasing to zero) and applying Lemma \ref{lem:convergence} with $H_n:=h_n^{-1}$, $\eta_n:=h_n^{-1}(\hat\mu_{[0,u_0+h_n],\cdot}-\hat\mu_{[0,u_0],\cdot})$ and $v:=\nu-\hat\nu_{[0,u_0],\cdot}$ completes the proof.
	\end{proof}
\section{Doob-like decomposition of the shadow couplings}\label{sec:martingale_points}
Fix $\mu,\nu\in\sP$ with $\mu\leq_{cd}\nu$. Let $\hat\mu\in\Pi(\lambda,\mu)$. The goal of this section is to determine the `martingale points' of an arbitrary lifted supermartingale shadow coupling. 

Recall that if $\hat\pi\in\hat\Pi_S(\hat\mu,\nu)$, then we can represent the coupling $\hat\pi$ as $\hat\pi(du,dx,dy)=\hat\pi_{(u,x),\cdot}(dy)\hat\mu(du,dx)$ where $\int_\R yd\hat\pi_{(u,x),\cdot}(dy)\leq x$ for $\hat\mu$-a.e. $(u,x)\in[0,1]\times\R$. On the other hand, if $\hat\pi$ is a lifted supermartingale shadow coupling, then by Theorem \ref{thm:shadowCouplings} we have that $\hat\pi(du,dx,dy)=\pi_{\hat\mu_{u,\cdot},T(u)}(dx,dy)du$, where (for $\lambda$-a.e. $u\in[0,1]$) the hitting coupling $\pi_{\hat\mu_{u,\cdot},T(u)}$ is the unique element of $\Pi_S(\hat\mu_{u,\cdot},\hat\mu_{u,\cdot}D_{T(u)})$. Recall that $\hat\nu_{u,\cdot}:=\hat\mu_{u,\cdot}D_{T(u)}$ is the (right) derivative of the curve $(S^\nu(\hat\mu_{[0,u],\cdot}))_{u\in[0,1]}$. Then $\hat\pi_{(u,x),\cdot}$ corresponds to the disintegration of the hitting coupling $\pi_{\hat\mu_{u,\cdot},T(u)}$ with respect to the first marginal $\hat\mu_{u,\cdot}$. It follows that (for $\lambda$-a.e. $u\in[0,1]$) $\pi_{\hat\mu_{u,\cdot},T(u)}\in\Pi_M(\hat\mu_{u,\cdot},\hat\nu_{u,\cdot})$ if and only if $\int_\R y\hat\pi_{(u,x),\cdot}(dy)=x$ for $\hat\mu_{u,\cdot}$-a.e. $x\in\R$. But, since $\pi_{\hat\mu_{u,\cdot},T(u)}$ is a supermartingale coupling, it is a martingale coupling whenever the means of its marginal distributions are equal, i.e., whenever $\overline{\hat\mu_{u,\cdot}}=\overline{\hat\nu_{u,\cdot}}$.

In order to identify the points $u\in[0,1]$ for which $\overline{\hat\mu_{u,\cdot}}=\overline{\hat\nu_{u,\cdot}}$ consider $c:[0,1]\to[0,\overline\mu-\overline\nu]$ defined by 
\begin{equation}\label{eq:defnC}
c(u)=c_{  \hat\mu_{[0,u],\cdot}, \nu },\quad u\in[0,1],
\end{equation}
 where $c_{\eta,\chi}$ is as in \eqref{eq:cDef1}.
\begin{lemma}\label{lem:c_u}
Let $c:[0,1]\to[0,\infty)$ be given by \eqref{eq:defnC}. Then $c(\cdot)$ is continuous and non-decreasing.
\end{lemma}
\begin{proof}
Fix $u,v\in[0,1]$ with $u\leq v$. By definition of $c_{\eta,\chi}$ (see \eqref{eq:cDef1}) we have that
$$
c(u)=\sup_{k\in\R}\{C_{\hat\mu_{[0,u],\cdot}}(k)-C_\nu(k)\}\leq\sup_{k\in\R}\{C_{\hat\mu_{[0,v],\cdot}}(k)-C_\nu(k)\}=c(v),
$$
where we used that $\hat\mu_{[0,v],\cdot}\geq \hat\mu_{[0,u],\cdot}$ and thus $C_{\hat\mu_{[0,v],\cdot}}\geq C_{\hat\mu_{[0,u],\cdot}}$ on $\R$.

We now prove that $c(\cdot)$ is continuous. Let $\eta_{v-u}:=\hat\mu_{[0,v],\cdot}-\hat\mu_{[0,u],\cdot}$ and $\chi_{v-u}:=S^\nu(\hat\mu_{[0,v],\cdot})-S^\nu(\hat\mu_{[0,u],\cdot})=S^{\nu-S^\nu(\hat\mu_{[0,u],\cdot})}(\eta_{v-u})$. By Lemma \ref{lem:shadowDefn} we have that $c(v)-c(u)=\overline{\eta_{v-u}}-\overline{\chi_{v-u}}$. Note that $\eta_{v-u},\chi_{v-u}\in\sM$ and $\eta_{v-u}(\R)=\chi_{v-u}(\R)=v-u$. Hence both $\eta_{v-u}$ and $\chi_{v-u}$ weakly converge to the zero measure, when either $v\downarrow u$ or $u\uparrow v$. Hence to conclude that $\lim_{v\downarrow u}[c(v)-c(u)]=0$ and $\lim_{u\uparrow v}[c(v)-c(u)]=0$ it is enough to show that the first moments of $\eta_{v-u}$ and $\chi_{v-u}$ converge to zero when $v\downarrow u$ or $u\uparrow v$, respectively. But this follows by observing that $\eta_{v-u}\leq\mu,\chi_{v-u}\leq\nu$ and both $\mu$ and $\nu$ are integrable (indeed, one can adapt the arguments of the last paragraph of the proof of Lemma \ref{lem:shadowConvergence}).
\end{proof}

We are now ready to present the main result of this section. Proposition \ref{prop:martingale_points} shows how given the initial data (i.e., $\mu\leq_{cd}\nu$ and a lift $\hat\mu\in\Pi(\lambda,\mu)$) one can immediately identify the `martingale points' of the corresponding lifted supermartingale shadow coupling from the graph of function $c(\cdot)$.
\begin{proposition}\label{prop:martingale_points}
Fix $\mu,\nu\in\sP$ with $\mu\leq_{cd}\nu$ and consider a lift $\hat\mu\in\Pi(\lambda,\nu)$. Define
$$\hat M:=\{u\in[0,1]:c'(u)\textrm{ exists and }c'(0)=0\}.$$

Then, for any  lifted supermartingale shadow coupling $\hat\pi\in\hat\Pi_S(\hat\mu,\nu)$ (as in Theorem \ref{thm:shadowCouplings}), $\hat M$ is a $\lambda$-a.s. unique (Borel) set for which $\hat\pi\lvert_{\hat M\times\R\times\R}$ is a martingale.
\end{proposition}
\begin{proof}
By the definition of $c(\cdot)$ and Lemma \ref{lem:shadowDefn} we have that $c(u)=\overline{\hat\mu_{[0,u],\cdot}}-\overline{S^\nu(\hat\mu_{[0,u],\cdot})}$, $u\in[0,1]$, and therefore, for each $h>0$,
$$
\frac{c(u+h)-c(u)}{h}=\frac{\overline{\hat\mu_{[0,u+h],\cdot}-\hat\mu_{[0,u],\cdot}}}{h}-\frac{\overline{S^{\nu-\hat\mu_{[0,u],\cdot}}(\hat\mu_{[0,u+h],\cdot}-\hat\mu_{[0,u],\cdot})}}{h}.
$$

Now let $L\subseteq[0,1]$ be a set for which $\hat\mu_{u,\cdot}$ and $\hat\nu_{u,\cdot}$ exist. Recall that $\lambda(L)=1$. Then
$$
\overline{\hat\mu_{u,\cdot}}-\overline{\hat\nu_{u,\cdot}}=\lim_{h\downarrow 0}\frac{c(u+h)-c(u)}{h}=c'(u+),
$$
i.e., the right derivative of $c(\cdot)$ at $u$ exists. But $c(\cdot)$ is non-decreasing, and therefore differentiable almost everywhere on $[0,1]$. It follows that $\overline{\hat\mu_{u,\cdot}}-\overline{\hat\nu_{u,\cdot}}=c'(u)$ for all $u\in\hat L:=L\setminus N_c$ where $N_c:=\{u\in[0,1]:c'(u)\textrm{ does not exist}\}$ is a $\lambda$-null set. Hence, if $u\in\{v\in\hat L:c'(v)=0\}$, then $\Pi_S(\hat\mu_{u,\cdot},\hat\nu_{u,\cdot})=\Pi_M(\hat\mu_{u,\cdot},\hat\nu_{u,\cdot})$ is a singleton with a unique element $\pi_{\hat\mu_{u,\cdot},T(u)}$, and it follows that $\hat\pi$ is a martingale on $\hat M\times\R\times\R$. The ($\lambda$-a.s.) uniqueness of $\hat M$ is straightforward.
\end{proof}
\section{The geometric construction of $\pi_D$}
\label{sec:geometric}
We fix $\mu,\nu\in\sP$ with $\mu\leq_{cd}\nu$ and $\overline\nu<\overline\mu$ throughout this section. Our goal here is to give an explicit construction of the decreasing supermartingale coupling introduced by Nutz and Stebegg \cite{NutzStebegg:18}. We begin by introducing monotonicity properties of the support of this coupling.
\begin{defn}\label{defn:support_mon}
Let $(\Gamma, M)\in\sB(\R^2)\times\sB(\R)$. We say
\begin{enumerate}
\item $\Gamma$ is second-order right-monotone if for all $(x,y_1),(x,y_2),(x',y')\in\Gamma$ with $x'<x$ we have that $y'\notin(y_1,y_2)$;
\item $(\Gamma,M)$ is first-order left-monotone if for all $(x_1,y_1),(x_2,y_2)\in\Gamma$ with $x_1<x_2$ and $x_2\notin M$ we have that $y_1\leq y_2$.
\end{enumerate}
\end{defn}

The following defines and characterizes the decreasing supermartingale coupling; see Nutz and Stebegg \cite[Theorems 1.1, 1.2 and 1.3]{NutzStebegg:18}.
\begin{defn}\label{defn:decreasing_coupling}
The decreasing supermartingale coupling, denoted by $\pi^D$, is the unique element of $\Pi_S(\mu,\nu)$ which satisfies any, and then all of the following
\begin{enumerate}
\item for each $x\in\R$, $\pi^D$ transports $\mu\lvert_{[x,\infty)}$ to the shadow $S^\nu(\mu\lvert_{[x,\infty)})$;
\item for all Borel $f:\R^2\to\R$ such that $f(x_2,\cdot)-f(x_1,\cdot)$ is strictly decreasing and strictly convex for all $x_1<x_2$, and $\lvert f(x,y)\lvert\leq a(x)+b(y)$ for all $x,y\in\R$ and some $\mu,\nu$ integrable functions $a,b:\R\to\R$, respectively, we have that
 $$\int fd\pi^D=\inf_{\pi\in\Pi_S(\mu,\nu)}\int fd\pi;$$
\item there exists first-order left-monotone and second-order right-monotone $(\Gamma, M)\in\sB(\R^2)\times\sB(\R)$ such that $\pi^D$ is concentrated on $\Gamma$ and $\pi^D\lvert_{M\times\R}$ is a martingale.
\end{enumerate}
\end{defn}
Note that $\pi^D$ is obtained by working from right to left and using the shadow measure (see the first characterization of Definition \ref{defn:decreasing_coupling}). In terms of lifted measures, this corresponds to taking $\hat\mu\in\hat\Pi(\lambda,\mu)$ to be the decreasing quantile lift of $\mu$. Then by applying Theorem \ref{thm:shadowCouplings} we obtain existence of a coupling $\hat\pi\in\hat\Pi_S(\hat\mu,\nu)$, for which $\int^1_0\hat\pi du=\pi^D$ (see Lemma \ref{lem:integratedPi}). Hence a construction of $\hat\pi$ leads to an explicit construction of $\pi^D$.

Let $G=G_\mu$ be the right continuous quantile function of $\mu$. Let $\hat\mu^D\in\Pi(\lambda,\mu)$ be the decreasing quantile lift of $\mu$, so that $\hat\mu(du,dx)=du\delta_{G(1-u)}(dx)$, or equivalently, $\hat\mu^D_{u,\cdot}=\delta_{G(1-u)}$, $u\in[0,1]$. (Note that we could redefine $\hat\mu_{u,\cdot}^D$ on a $\lambda$-null set, and thus in fact we could use any version of a generalized inverse of $F_\mu$ to represent $\hat\mu^D$).

For each $u\in[0,1]$, define $\mu_u\in\sM$ by
\begin{equation}\label{eq:mu_u}\mu_{u}(A)= \mu\Big(A \cap \big(G(1-u),\infty\big)\Big)+\Bigg(u-\mu\Big(\big(G(1-u),\infty\big)\Big)\Bigg)\delta_{G(1-u)}(A),\quad \textrm{for all Borel }A\subseteq\R.\end{equation}
Then $\mu_u=\int^u_0\hat\mu^D_{v,\cdot}dv=\hat\mu^D_{[0,u],\cdot}$.

\begin{lemma}\label{lem:integratedPi}
	Let $\mu,\nu\in\sP$ with $\mu\leq_{cd}\nu$. Let $\pi^D$ be the decreasing supermartingale coupling of $\mu$ and $\nu$. Let $\hat\pi^D\in \hat\Pi_S(\hat\mu^D,\nu)$ be the unique lifted shadow coupling of $\hat\mu^D$ and $\nu$ (as in Theorem \ref{thm:shadowCouplings}). Then $\int^1_0\hat\pi^D du=\pi^D$.
\end{lemma}
\begin{proof}
By Theorem  \ref{thm:shadowCouplings}, $\hat\pi^D$ is a unique measure that, for each $u\in[0,1]$, transports $\mu_u$ to $S^\nu(\mu_u)$. On the other hand, fix $x\in \sI_\mu$ and let $u_x\in[0,1]$ be given by $u_x:=\sup\{u\in[0,1]: x\leq G(1-u)\}$. It follows that $\mu_{u_x}=\mu\lvert_{[x,\infty)}$, and then for all Borel $B\subseteq \R$ we have that
$$
\int_{[0,1]\times[x,\infty)\times B}\hat\pi^D(du,dx,dy)=\int_{[0,u_x]\times[x,\infty)\times B}\hat\pi^D(du,dx,dy)=S^\nu(\mu_{u_x})(B)=S^\nu(\mu\lvert_{[x,\infty)})(B),
$$
which shows that $\int^1_0\hat\pi^D du=\pi^D$.
\end{proof}

We now provide an explicit construction of $\hat\pi^D\in \hat\Pi_S(\hat\mu^D,\nu)$.

Recall the definition of $D(k)=P_\nu(k)-P_\mu(k)$, $k\in\R$, and that $D\geq0$ on $\R$. Note that, since $\overline\nu<\overline\mu$, $\lim_{k\to\infty}D(k)=\overline\mu-\overline\nu>0$. In what follows (and in the light of Section \ref{sec:irreducible}) we assume that $\{k\in\R:D(k)>0\}=(\ell_D,r_D)=(\ell_D,\infty)$ is an (open) interval, $\mu((\ell_D,\infty))=1$ and $\nu((\ell_D,\infty))+\nu(\{\ell_D\})=1$ with $\nu(\{\ell_D\})=0$ whenever $\ell_D=-\infty$.

Recall also the definition of the sub-differential $\partial h(x)$ of a convex function $h:\R \mapsto \R$ at $x$:
\[ \partial h (x) = \{ \phi \in\R: h(y) \geq h(x) + \phi(y-x) \mbox{ for all } y \in \R \}. \]
We extend this definition to non-convex functions $f$ so that the subdifferential of $f$ at $x$ is given by
\[ \partial f (x) = \{ \phi\in\R : f(y) \geq f(x) + \phi(y-x) \mbox{ for all } y \in \R \}. \]
If $h$ is convex then $\partial h$ is non-empty everywhere, but this is not the case for non-convex functions. Instead we have that $\partial f(x)$ is non-empty if and only if $f(x)=f^c(x)$ and then $\partial f^c(x) = \partial f(x)$.

For each $u\in[0,1]$, let $\mu_u\in\sM$ be defined as in \eqref{eq:mu_u}. We have $C_{\mu_u}(k)=C_\mu(k)$ for $k\geq G(1-u)$, while $C_{\mu_u}(k)\leq C_\mu(k)$ for $k< G(1-u)$. In particular,
 \begin{equation*}
 C_{\mu_u}(k)=C_\mu(k\vee G(1-u))-u(G(1-u)-k)^+,\quad k\in\R,
 \end{equation*}
and thus $C_{\mu_u}(\cdot)$ is linear on $(-\infty,G(1-u))$ and $(-u)\in\partial C_\mu(G(1-u))$, so that $C'_\mu(G(1-u)-)\leq -u\leq C'_\mu(G(1-u)+)$.

For each $u \in [0,1]$ define $\sE_u:\R\mapsto\R_+$ by $\sE_u=D+C_\mu-C_{\mu_u}$. Note that $\sE_u(k)=D(k)$ for $k\geq G(1-u)$. Since $C_\mu-C_{\mu_u}$ is non-negative on $\R$, we have that $\sE_u(k)\geq D(k)$ for $k< G(1-u)$. Moreover, since $C_{\mu_u}$ is linear on $(-\infty,G(1-u))$, $\sE_u$ is convex on $(-\infty,G(1-u))$. It is also easy to see that both $u\mapsto\sE_u(k)$ (for a fixed $k\in\R$) and $k\mapsto\sE_u(k)-D(k)$ (for a fixed $u\in[0,1]$) are non-increasing.

Recall that by Lemma \ref{lem:shadowDefn}
 \begin{equation*}
 P_{S^\nu(\mu_u)}(k)=P_\nu(k)-(P_\nu-P_{\mu_u})^c(k),\quad k\in\R.
 \end{equation*}
Next lemma shows that we can also identify $S^\nu(\mu_u)$ by considering the convex hull of $\sE_u$.
\begin{lemma}\label{lem:C_shadow}
	Let $\mu,\nu\in\sP$ with $\mu\leq_{cx}\nu$. Consider $(\mu_u)_{u\in[0,1]}$ where $\mu_u$ is defined as in \eqref{eq:mu_u}. Let $c:[0,1]\to\R$ be as in \eqref{eq:defnC}, i.e., $c(u)=c_{\mu_u,\nu}$, $u\in[0,1]$. 
Then, for each $u\in[0,1]$,
	\begin{align*}
	C_{S^\nu({\mu_u})}(k)&=C_\nu(k)-\sE^c_u(k)+(\overline\mu-\overline\nu)-(\overline{\mu_u}-\overline{S^\nu(\mu_u)})\\&=C_\nu(k)-\sE^c_u(k)+c(1)-c(u),\quad k\in\R.
	\end{align*}
	In particular, $S^\nu(\mu_u)$ corresponds to second (distributional) derivative of $(C_\nu-\sE^c_u)$.
\end{lemma}
\begin{proof}
	Using the Put-Call parity (twice) and  Lemma \ref{lem:shadowDefn} we have that
\begin{align}
C_{S^\nu(\mu_u)}(k)&= P_{S^\nu(\mu_u)}(k)+\overline{S^\nu(\mu_u)}-S^\nu(\mu_u)(\R)k\nonumber\\
&=(C_\nu(k)-\overline\nu+\nu(\R)k)-(P_\nu-P_{\mu_u})^c(k)+(\overline{S^\nu(\mu_u)}-S^\nu(\mu_u)(\R)k),\quad k\in\R.\label{eq:C_shadow}
\end{align}
 On the other hand, by the Put-Call parity and definition of $\sE_u$ we also have that
 $$P_\nu(k)-P_{\mu_u}(k)=\sE_u(k)-(\overline\mu-\mu(\R)k)+(\overline{\mu_u}-uk),\quad k\in\R.$$
 Then by Beiglb\"{o}ck et al. \cite[Lemma 3]{BHN:20} and linearity of $k\mapsto\left\{(\overline{\mu_u}-uk)-(\overline\mu-\mu(\R)k)\right\}$, $(P_\nu-P_{\mu_u})^c(k)=\sE^c_u(k)-(\overline\mu-\mu(\R)k)+(\overline{\mu_u}-uk)$, which together with \eqref{eq:C_shadow} and definition of $c(\cdot)$ proves the claim (here we used that $\mu(\R)=\nu(\R)$ and $u=\mu_u(\R)=S^\nu(\mu_u)(\R)$).
\end{proof}

Since $\hat\mu^D_{u,\cdot}=\delta_{G(1-u)}$, $u\in[0,1]$, the hitting coupling $\pi_{\hat\mu_{u,\cdot},T(u)}$, where $T(u)=supp(\nu-S^\nu(\hat\mu_{[0,u],\cdot}))$, is in fact a product measure of $\delta_{G(1-u)}$ and the hitting projection $(\delta_{G(1-u)}D_{T(u)})$. In particular,
\begin{align}\label{eq:hittingPiD}&\pi_{\hat\mu_{u,\cdot},T(u)}(dx,dy)\nonumber\\=&\delta_{G(1-u)}(dx)(\delta_{G(1-u)}D_{T(u)})(dy)\nonumber\\=&\begin{cases}
\delta_{G(1-u)}(dx)\delta_{G(1-u)}(dy)& \textrm{if }G(1-u)\in T(u);\\
\delta_{G(1-u)}(dx)\left[\frac{s(u)-G(1-u)}{s(u)-r(u)}\delta_{r(u)}(dy)+\frac{G(1-u)-r(u)}{s(u)-r(u)}\delta_{s(u)}(dy)\right]&\textrm{otherwise};
\end{cases}\\
=&\delta_{G(1-u)}(dx)\chi_{r(u),G(1-u),s(u)}(dy),\nonumber
\end{align}
where
\begin{align}
r(u)&=\sup\{k\in T(u)\cap(-\infty, G(1-u)]\},\label{eq:f}\\
s(u)&=\inf\{k\in T(u)\cap[G(1-u),\infty)\}.\label{eq:g}
\end{align}
Our goal is using Lemma \ref{lem:C_shadow} to identify the versions of $r$ and $s$ from the graphs of $\sE^c_u$, $u\in[0,1]$.

\begin{defn}
	\label{def:slope}
	$\phi:(0,1) \mapsto \R$ is given by $\phi(u) = \sup \{ \psi: \psi \in \partial \sE_u^c (G(1-u)) \}=(\sE_u)'(G(1-u)+)$.
\end{defn}

Recall the definition of $ L^f_{a,b}$ for any $f:\R\mapsto\R$ (see \eqref{eq:L1}), so that (in the case $a<b$) $L^f_{a,b}$ is the line passing through $(a,f(a))$ and $(b, f(b))$. Define also $L^{f,\psi}_a$ by $L^{f,\psi}_a(y) = f(a) + \psi(y-a)$ so that $L^{f,\psi}_a$ is the line passing through $(a,f(a))$ with slope $\psi$. (Note that, in the case $a=b$, $L^f_{a,a}=L^{f,0}_a$.)

Define $R,S:[0,1]\mapsto\R\cup\{-\infty,\infty\}$ by
\begin{align}  
R(u)&:=\inf\{k\in\R:k\leq G(1-u), \sE^c_u(k)= {L}^{\sE^c_u ,\phi(u)}_{G(1-u)}(k)\},\quad u\in[0,1],  \label{eq:R}\\
S(u)&:=  \sup\{k\in\R:k \geq G(1-u), \sE^c_u(k) = {L}^{\sE^c_u ,\phi(u)}_{G(1-u)}(k) \},\quad u\in[0,1].\label{eq:S}
\end{align}
See Figure \ref{fig:RS}.

We first establish global monotonicity properties of $(R,S)$. See Figure \ref{fig:RSglobal}.
\begin{proposition}\label{prop:RSmonotonicity}
Let $R$ and $S$ be as in \eqref{eq:R} and \eqref{eq:S}, respectively. Then the pair $(R,S)$ is first-order left-monotone and second order right-monotone with respect to $G$ on $[0,1]$:
	\begin{enumerate}
		\item $R(u)\leq G(1-u)\leq S(u)$ for all $u\in[0,1]$;
		\item $R$ is decreasing on $[0,1]$;
		\item For $u,v\in[0,1]$ with $u<v$, $S(v)\notin(R(u),S(u))$;
		\item  For $u,v\in[0,1]$ with $u<v$, $\sup\{k\in T(v)\}\leq\inf\{k\in T(u)\}$ on $\{u\in[0,1]:S(u)=\infty\}$.
	\end{enumerate}
\end{proposition}
\begin{figure}[H]
	\centering

	\caption{Plot of locations of $R(u)$, $G(1-u)$, $S(u)$, $R(v)$ and $G(1-v)$, for $u<v$, and in the case where $-\infty<R(u)<G(1-u)<S(u)<\infty$ and $-\infty<R(v)<G(1-v)<S(v)=\infty$. The dashed curve represents $D$. Note that $\lim_{k\to\infty}D(k)=\overline\mu-\overline\nu>0$. The dotted curves correspond to the graphs of $\sE_{u}$ and $\sE_{v}$. Note that $D=\sE_u$ on $[G(1-u),\infty)$ (resp. $D=\sE_v$ on $[G(1-v),\infty)$), while $\sE_u$ (resp. $\sE_v$) is convex and $D\leq\sE_u$ (resp. $D\leq\sE_v$ ) on $(-\infty,G(1-u))$ (resp. $(-\infty,G(1-v))$). The solid curves below $\sE_{u}$ and $\sE_v$ represent $\sE^c_{u}$ and $\sE^c_{v}$, respectively. The convex hull $\sE^c_{u}$ (resp. $\sE^c_{v}$) is linear on $(R(u),S(u))$ (resp. $(R(u),\infty)$).}
	\label{fig:RS}
\end{figure}
\begin{proof}
Fix $u,v\in[0,1]$ with $u<v$.
	
	(1) This is immediate from the definitions of $R$ and $S$.
	
	(2)-(3) If $R(u)=G(1-u)=S(u)$, then $(R(u),S(u))=\emptyset$ and $S(v)\notin(R(u),S(u))$ by default. Also, $R(v)\leq G(1-v)\leq G(1-u)=R(u)$, as required.
	
	Hence suppose that $R(u)<S(u)$. We have that $\sE^c_u$ is linear on $(R(u),S(u))$ and therefore $(\nu-S^\nu(\mu_u))$ does not charge $(R(u),S(u))$. By the associativity of the shadow measure (see \eqref{eq:shadowAssoc} in Lemma \ref{lem:shadowDefn}), $\nu-S^\nu(\mu_v)=\nu-S^\nu(\mu_u)-S^{\nu-S^\nu(\mu_u)}(\mu_v-\mu_u)$ and therefore  $(\nu-S^\nu(\mu_v))$ does not charge  $(R(u),S(u))$ as well. It follows that $\sE^c_v$ is linear on $(R(u),S(u))$.
	
	Now suppose $S(v)\in (R(u),S(u))$. By the definition of $S$ and convexity of $\sE^c_v$ we have that $\sE^c_v>L^{\sE^c_v,\phi(v)}_{G(1-v)}$ on $(S(v),\infty)$. It follows that the second derivative of $\sE^c_v$ corresponds to a measure $\eta_v\in\sM$ with $\eta_v([S(v),S(u)))>0$. But by Lemma \ref{lem:C_shadow} we have that $(\nu-S^\nu(\mu_v))=\eta_v$; a contradiction since $(\nu-S^\nu(\mu_v))$ does not charge  $(R(u),S(u))$ and thus also $[S(v),S(u))$.
	
	We now show that (in the case $R(u)<S(u)$) $R(v)\leq R(u)$. Suppose not, so that $R(u)<R(v)$. If $R(u)<R(v)<S(u)$, then (similarly as in the case for $S$) we must have that $(\nu-S^\nu(\mu_v))((R(u),R(v)])>0$ which contradicts the fact that $(\nu-S^\nu(\mu_v))$ does not charge  $(R(u),S(u))$ and thus also $(R(u),R(v)]$. Hence we assume that $R(v)\geq S(u)$. Then
	\begin{equation}\label{eq:proofRGS}
	G(1-u)\leq S(u)\leq R(v)\leq G(1-v).
	\end{equation} Since $G$ is non-decreasing and $u<v$, we have a contradiction if at least one inequality in \eqref{eq:proofRGS} is strict. Therefore $R(u)<G(1-u)= S(u)= R(v)= G(1-v)$.
	
	Now note that, since $R(v)=G(1-v)$, we must have that $\sE^c_v=\sE^c_u$ on $[G(1-v)=G(1-u),\infty)$ and therefore $\phi(v)=\phi(u)$. It follows that $L^{\sE^c_v,\phi(v)}_{G(1-v)}=L^{\sE^c_u,\phi(u)}_{G(1-u)}$. By convexity of $\sE^c_v$ we have that $\sE^c_v\geq L^{\sE^c_v,\phi(v)}_{G(1-v)}$ on $(-\infty,G(1-v)]$. On the other hand, $\sE_u\geq\sE_v$ on $\R$, and therefore $\sE^c_u\geq\sE^c_v$ on $\R$. It follows that
	$$
	\sE^c_u=L^{\sE^c_u,\phi(u)}_{G(1-u)}=L^{\sE^c_v,\phi(v)}_{G(1-v)}\leq \sE^c_v,\quad\textrm{on } [R(u),G(1-v)],
	$$
	and therefore $\sE^c_u=L^{\sE^c_v,\phi(v)}_{G(1-v)}=\sE^c_v$ on $[R(u),G(1-v)]$. But then $R(v)>R(u)\geq \inf\{k\leq G(1-v):\sE^c_v(k)= L^{\sE^c_v,\phi(v)}_{G(1-v)}(k)\}=R(v)$, a contradiction. We conclude that $R(v)\leq R(u)$.
	
	(4) Finally suppose that $R(u)\leq G(1-u)<S(u)=\infty$. Then $\sE^c_u$ is linear on $(R(u),\infty)$ and $(\nu-S^\nu(\mu_u))$ does not charge $(R(u),\infty)$. It follows that $\sup\{k\in T(u)\}=R(u)$. By the associativity of the shadow measure we have that $(\nu-S^\nu(\mu_v))$ does not charge $(R(u),\infty)$ as well, and therefore $\sup\{k\in T(v)\}\leq\sup\{k\in T(u)\}=R(u)$ as required.
\end{proof}

The following is the main result of this section.
\begin{theorem}\label{thm:construction}
	Let $\mu,\nu\in\sP$ with $\mu\leq_{cd}\nu$. Let $(R,S)$ be given by \eqref{eq:R} and \eqref{eq:S}, and define $\hat\pi^{R,S}\in\sM^3$ by
	$$
	\hat\pi^{R,S}(du,dx,dy)=du\delta_{G(1-u)}(dx)\chi_{R(u),G(1-u),S(u)}(dy).
	$$
	Then $\hat\pi^{R,S}=\hat\pi^D$, so that $\int^1_0\hat\pi^{R,S} du$ is the decreasing supermartingale coupling $\pi^D$.

\end{theorem}
\begin{proof}
	By \eqref{eq:hittingPiD} it is enough to show that $\chi_{R(u),G(1-u),S(u)}=\chi_{r(u),G(1-u),s(u)}$.
	
We claim that 

\begin{align*}\{u\in[0,1]:G(1-u)\notin T(u)\}&=\{u\in[0,1]:r(u)<G(1-u)<s(u)\}\\&=\{u\in[0,1]:R(u)<G(1-u)<S(u)\}\\
&=\{u\in[0,1]:G(1-u)\notin T(u), R(u)=r(u), S(u)=s(u)\}.\end{align*}
Note that the first equality is an immediate consequence of the definitions of $T(u),r(u),s(u)$. We now simultaneously establish the second and third equalities.

First, let $u\in[0,1]$ be such that $r(u)<G(1-u)<s(u)$. Then $(\nu-S^\nu(\mu_u))$ does not charge $(r(u),s(u))$. By Lemma \ref{lem:C_shadow} we then have that $\sE^c_u$ is linear on $(r(u),s(u))$. But by the definitions of $r$ and $s$, $(r(u),s(u))$ is the largest open interval $I\ni G(1-u)$ with $(\nu-S^\nu(\mu_u))(I)=0$. Consequently, $(r(u),s(u))$ is also the largest open interval $\tilde I\ni G(1-u)$ such that $\sE^c_u$ is linear on $\tilde I$. It follows that $R(u)=r(u)$ and $S(u)=s(u)$.
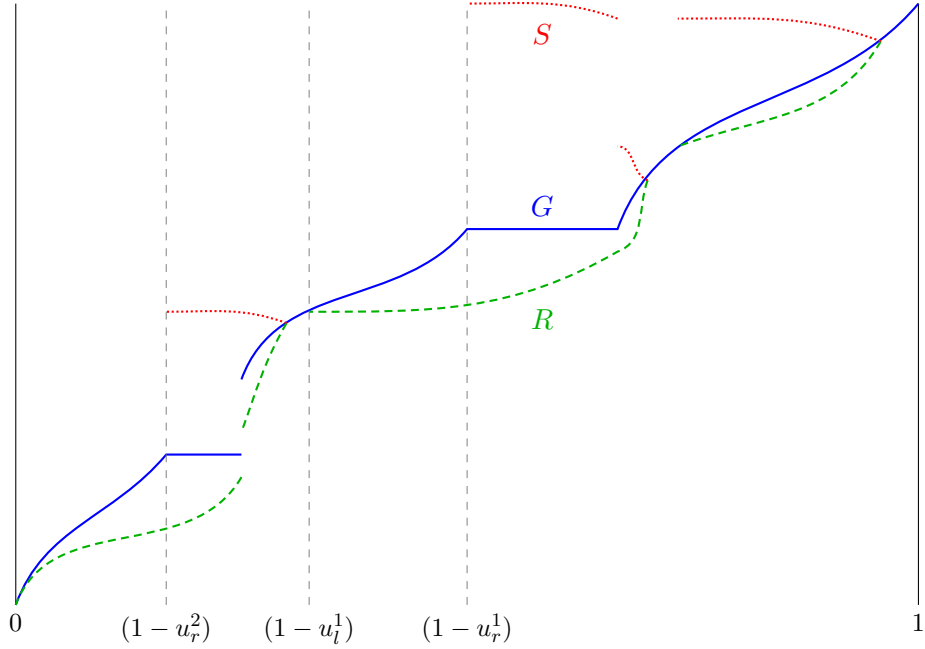
\begin{figure}[H]
	\centering
	\begin{tikzpicture}[scale=1,
	declare function={	
		k1=2.1;
		k2=1;
		a(\x)=((k1-\x)*(\x<k1))-k1-1;
		b(\x)=((k2-\x)*(\x<k2))-k1-1;
		x1=-6.8;
		x2=-6.3;
		z1=-5.2;
		z2=-4.7;
	}]

	\draw[ black] (-9,-3)--(-9,5);
	\draw[black] (3,-3)--(3,5);
	

	\draw[blue,thick, name path=g1] (-9,-3) to[out=70, in=230] (-7,-1) -- (-6,-1) ;
	\draw[blue,thick, name path=g2] (-6,0)  to[out=70, in=230] (-3,2) -- (-1,2) to[out=70, in=230] (3,5);
	
	\draw[black!30!green,thick, densely dashed, name path=r2] (2.5,4.5) to[out=240, in=20] (-.2,3.1) ;
	
	\draw[red,thick, densely dotted, name path=r3] (2.5,4.5) to[out=160, in=0] (-.2,4.8) ;
	
	\draw[black!30!green,thick, densely dashed, name path=r2] (-.6,2.65) to[out=250, in=20] (-1,1.7) to[out=210, in=0] (-5.1,0.9);
	
	\draw[red,thick, densely dotted, name path=r3] (-.6,2.65) to[out=160, in=0] (-1,3.1) ;
	
	\draw[red,thick, densely dotted, name path=r3] (-1,4.8) to[out=160, in=0] (-3,5) ;
	
	\draw[black!30!green,thick, densely dashed, name path=r2] (-5.4,0.75) to[out=240, in=70] (-6,-0.7) ;
	
	\draw[red,thick, densely dotted, name path=r3] (-5.4,0.75) to[out=160, in=0] (-7,0.9) ;
	
	\draw[black!30!green,thick, densely dashed, name path=r2] (-6,-1.3) to[out=240, in=70] (-9,-3) ;
	
	\draw[gray, very thin, dashed] (-7,-3) -- (-7,5);
	\draw[gray, very thin, dashed] (-5.1,-3) -- (-5.1,5);
	\draw[gray, very thin, dashed] (-3,-3) -- (-3,5);

	\node[red] at (-2,4.6) {$S$};
	\node[blue] at (-2,2.3) {$G$};
	\node[black!30!green] at (-2,0.8) {$R$};
	
	\node[below,scale=0.9] at (-7,-3) {$(1-u^2_r)$};
	\node[below,scale=0.9] at (-5.1,-3) {$(1-u^1_l)$};
	\node[below,scale=0.9] at (-3,-3) {$(1-u^1_r)$};
	\node[below,scale=0.9] at (-9,-3) {$0$};
	\node[below,scale=0.9] at (3,-3) {$1$};
	
	\end{tikzpicture}
	\caption{Sketch of the supporting functions of the lifted decreasing supermartingale coupling $\hat\pi^D$: $R$ (dashed), $G$ (solid) and $S$ (dotted). In the figure, $I_1:=(1-u^1_l,1-u^1_r)$ and $I_2:=(0,1-u^2_r)$ are the strict supermartingale regions of $\hat\pi^D$. In particular, for $u\in[0,1]$ with $(1-u)\in (I_1\cup I_2)$, the mass at $G(1-u)$ is mapped to $R(u)$ only, and thus $\hat\pi^D$ resembles the deterministic (lifted) Hoeffding-Fr\'echet coupling $\hat\pi^{HF}$. On the other hand, for $u\in[0,1]$ with $(1-u)\notin(I_1\cup I_2)$, $\hat\pi^D$ maps the mass at $G(1-u)$ to two points $\{R(u),S(u)\}$, and thus resembles the (lifted) right-curtain martingale coupling $\hat\pi^{rc}$.}
	\label{fig:RSglobal}
\end{figure}
Conversely, suppose $u\in[0,1]$ is such that $R(u)<G(1-u)<S(u)$. Then $\sE^c_u$ is linear on $(R(u),S(u))$ and by Lemma \ref{lem:C_shadow} we have that $(\nu-S^\nu(\mu_u))$ does not charge $(R(u),S(u))$. But $(R(u),S(u))$ is the largest open interval $\tilde I\ni G(1-u)$ such that $\sE^c_u$ is linear on $\tilde I$. Consequently, $(R(u),S(u))$ is the largest open interval $I\ni G(1-u)$ with $(\nu-S^\nu(\mu_u))(I)=0$. It follows that $R(u)=r(u)$ and $S(u)=s(u)$. 

We conclude that $\chi_{R(u),G(1-u),S(u)}=\chi_{r(u),G(1-u),s(u)}$ whenever $G(1-u)\notin T(u)$, $u\in[0,1]$.

Furthermore,
\begin{align*}
&\{u\in[0,1]:G(1-u)\in T(u)\}\\=&[0,1]\setminus\{u\in[0,1]:G(1-u)\notin T(u)\}\\=&\{u\in[0,1]:R(u)=G(1-u)<S(u)\}\cup\{u\in[0,1]:R(u)<G(1-u)=S(u)\}\\&\quad\quad\cup\{u\in[0,1]:R(u)=G(1-u)=S(u)\}.
\end{align*}
But if $R(u)=G(1-u)$ or $S(u)=G(1-u)$ (or both) then $$\chi_{R(u),G(1-u),S(u)}=\delta_{G(1-u)}=\chi_{r(u),G(1-u),s(u)},$$
where the second equality follows from the fact that $r(u)=G(1-u)=s(u)$ whenever $G(1-u)\in T(u)$, $u\in[0,1]$.

Combining both cases we conclude that $\chi_{R(u),G(1-u),S(u)}=\chi_{r(u),G(1-u),s(u)}$.
\end{proof}

\begin{remark}\label{rem:mainTheorem} The proof of Theorem \ref{thm:construction} can also be obtained by using Proposition \ref{prop:RSmonotonicity} together with the third characterization of $\pi^D$ in Definition \ref{defn:decreasing_coupling}.
	
The first two properties of $R$ and $S$ in Proposition \ref{prop:RSmonotonicity} translate to the second-order right-monotonicity of the coupling $\int^1_0\hat\pi^{R,S}du$ with respect to $\Gamma:=\bigcup_{u\in[0,1]}\Big\{(G(1-u),R(u))\cup(G(1-u),G(1-u))\cup(G(1-u),S(u))\Big\}\in\sB(\R^2)$ (in the sense of Definition \ref{defn:support_mon}).

Furthermore, for $u\in[0,1]$, $\chi_{R(u),G(1-u),S(u)}$ is a strict supermartingale kernel (i.e., $G(1-u)>\overline{\chi_{R(u),G(1-u),S(u)}}$) if and only if $R(u)<G(1-u)<S(u)=\infty$. Hence $\tilde M:=\{u\in[0,1]: S(u)<\infty\}\cup\{u\in[0,1]:R(u)=G(1-u)<S(u)=\infty\}$ is such that $\hat\pi^{R,S}\lvert_{\tilde{M}\times\R\times\R}$ is a martingale. Furthermore, it is easy to see that $\lambda(\hat M\cap\tilde M)=\lambda (\hat M)=\lambda (\tilde M)$, where $\hat M$ is as in Proposition \ref{prop:martingale_points}. This together with property (3) of Proposition \ref{prop:RSmonotonicity} imply the first-order left-monotonicity of $\int^1_0\hat\pi^{R,S}du$ with respect to $(\Gamma, M:=\{G(1-u):u\in \tilde M\})$ (in the sense of Definition \ref{defn:support_mon}).

Using the characterization of $\pi^D$ (in terms of the monotonicity of its support) one then establishes that $\int^1_0\hat\pi^{R,S}du=\pi^D$.
\end{remark}

\appendix
\section{Up and Down measures}
Consider $\chi,\chi'\in\sM$ with $\chi(\R)=\chi'(\R)$, and define Up$(\chi,\chi')$, Down$(\chi,\chi')\in\sM$ as in \eqref{eq:Up} and \eqref{eq:Down}, respectively.
\begin{lemma}\label{lem:UpDownSto}
For $\chi,\chi'\in\sM(\R)$ with $\chi(\R)=\chi'(\R)$ the following holds
$$
\textrm{\normalfont Down}(\chi,\chi')\leq_{sto}\chi \leq_{sto}\textrm{\normalfont Up}(\chi,\chi')\quad\textrm{and}\quad\textrm{\normalfont Down}(\chi,\chi')\leq_{sto}\chi'\leq_{sto} \textrm{\normalfont Up}(\chi,\chi').
$$
\end{lemma}
\begin{proof}
We only prove that Down$(\chi,\chi'))\leq_{sto}\chi$. The other relations use similar arguments. It is enough to show that $F_{\textrm{Down}(\chi,\chi')}\geq F_\chi$ everywhere: for each $x\in\R$,
\begin{align*}
F_{\textrm{Down}(\chi,\chi')}(x)&=\lambda\lvert_{[0,\chi(\R)=\chi'(\R)]}(\{u\in\R:(G_\chi(u)\wedge G_{\chi'}(u))\leq x\})\\
&\geq\lambda\lvert_{[0,\chi(\R)]}(\{u\in\R:G_\chi(u)\leq x\})=F_\chi(x).
\end{align*}
\end{proof}
\begin{lemma}[{Juillet \cite[Lemma 2.25]{Juillet:16}}]\label{lem:UpDownW}
For $\chi,\chi'\in\sM$ with $\chi(\R)=\chi'(\R)$ the following holds
\begin{align*}
W(\chi,\chi')&=W(\chi,\textrm{\normalfont Down}(\chi,\chi'))+W(\chi',\textrm{\normalfont Down}(\chi,\chi'))\\
&=W(\chi,\textrm{\normalfont Up}(\chi,\chi'))+W(\chi',\textrm{\normalfont Up}(\chi,\chi')).
\end{align*}
\end{lemma}
\begin{lemma}\label{lem:UpDownPCD}
Consider $\chi,\chi'\in\sM$ with $\chi(\R)=\chi'(\R)$. Let $\eta\in\sM$ be such that $\eta\leq_{pcd}\chi$ and $\eta\leq_{pcd}\chi'$. Then $\eta\leq_{pcd}$ {\normalfont Down}$(\chi,\chi')$.
\end{lemma}
\begin{proof}
By Lemma \ref{lem:UpDownSto} we have that Down$(\chi,\chi')\leq_{sto}\chi$ and  Down$(\chi,\chi')\leq_{sto}\chi'$. Then, for any non-increasing $f:\R\to\R$, we have that $\left(\int_\R f d\chi\right)\vee\left(\int_\R f d\chi'\right)\leq\int_\R fd$Down$(\chi,\chi')$, where we used that $(-f)$ is non-decreasing. Clearly any  positive, convex and non-increasing $g:\R\to\R$ is also non-increasing, and thus the claim follows.
\end{proof}
\begin{lemma}\label{lem:Up1marginals}
Let $\eta,\eta',\chi\in\sM$ be such that $\eta(\R)=\eta'(\R)$, $\eta\leq_{pcd}\chi$ and $\eta'\leq_{pcd}\chi$. Then {\normalfont Up}$(\eta,\eta')\leq_{pcd} \chi$.
\end{lemma}
\begin{proof}
Let $f:\R\to\R$ be non-negative, convex and non-increasing. By Lemma \ref{lem:UpDownSto} we have that $\eta\leq_{sto}$Up$(\eta,\eta')$, from which we deduce that $\int_\R(-f)d\eta\leq\int_\R(-f)d\textrm{Up}(\eta,\eta')$, where we used that $(-f)$ is non-decreasing. It follows that $\int_\R fd\textrm{Up}(\eta,\eta')\leq\int_\R fd\eta\leq\int_\R fd\chi$.
\end{proof}
\section{Convex hull}
Let $f:\R\to\R$ be measurable, and denote by $f^c$ the convex hull of $f$ (i.e., the largest convex function below $f$). Note that we may have that $f^c=-\infty$ on $\R$. Furthermore, if a function $f=-\infty$ (or $f=\infty$) on $\R$ then we deem it to be both linear and convex, and then set $f^c=f$.

Fix $x,z\in\R$ with $x \leq z$, and define $L^f_{x,z}:\R\mapsto\R$ by
\begin{equation}\label{eq:L1}
L^f_{x,z}(y)=\begin{cases}
f(x)+\frac{f(z)-f(x)}{z-x}(y-x),&\textrm{ if }x<z,\\
f(x),&\textrm{ if }x=z.
\end{cases}
\end{equation}
Then (see Rockafellar \cite[Corollary 17.1.5]{Rockafellar:72}),
\begin{equation}\label{eq:chull}
f^c(y)=\inf_{x\leq y\leq z}L^f_{x,z}(y),\quad y\in\R.
\end{equation}
(Note that for \eqref{eq:chull}, the definition of $L^f_{x,z}$ outside $[x,z]$ is irrelevant and we could restrict the domain of $L^f_{x,z}$ to $[x,z]$.)
\begin{lemma}\label{lem:hullDerivative}
Let $f,g:\R\to\R$ be the differences of two convex functions. If $f'_+\leq g'_+$ then $(f^c)'_+\leq(g^c)'_+$ on $\R$.
\end{lemma}
\begin{proof}
	Note that for (a measurable) $h:\R\to\R$ and $a\in\R$ we have
	$$
	h^c(x)+a=(h+a)^c(x),\quad x\in\R,
	$$
	and therefore $(h^c)'_+(x)=(h^c+a)'_+(x)$ for each $x\in\R$.
	
	Fix $x\in\R$. Using the above observation, without loss of generality we can assume that $f(x)=g(x)$. Then since $f'_+\leq g'_+$ on $\R$,
	$$
	f\geq g\textrm{ on }(-\infty,x),\quad f\leq g\textrm{ on }(x,\infty)
	$$
	(this can be easily proved by using the absolute continuity of $f$ and $g$).
	
	We will prove the claim by contradiction; suppose that $(f^c)'_+(x)>(g^c)'_+(x)$. There are three cases.
	
	1. Suppose $f^c(x)=g^c(x)$. Then $g\geq f\geq f^c\geq l^{f^c}_x>l^{g^c}_x$ on $(x,\infty)$, where $z\to l^h_x(z):=h(x)+h'_+(x)(z-x)$. It follows that
	$$
	z\to\tilde g(z):=g^c(z)I_{\{z\leq x\}}+l^{f^c}_x(z)I_{\{z>x\}}
	$$ 
	is a convex minorant of $g$, and therefore $g^c\geq\tilde g$ on $\R$. But then, since $g^c(x)=\tilde{g}(x)$,
	$$
	(g^c)'_+(x)\geq\tilde{g}'_+(x)=(f^c)'_+(x)>(g^c)'_+(x),
	$$
	a contradiction.
	
	2. Suppose $g^c(x)>f^c(x)$. Then $f(x)=g(x)\geq g^c(x)>f^c(x)$, and we have that there exists an interval $I\subseteq\R$ with $x\in I^\circ$ such that $f^c$ is linear on $I$ (see, for example, Hobson and Norgilas \cite[Lemma 2.2]{HobsonNorgilas:21}), so that $f^c=l^{f^c}_x$ on $I$. Define $z\to\tilde l(z):=l^{g^c}_x(z)-g^c(x)+f^c(x)$, and note that $\tilde l(x)=l^{f^c}_x(x)$ and $l^{f^c}_x<\tilde l$ on $(-\infty, x)$.
	
	Suppose there exists $\underline x<x$ such that $f^c(\underline x)=\tilde{l}(\underline x)$. Then, since $f\geq g\geq g^c\geq l^{g^c}_x>\tilde l$ on $(-\infty,x)$,
	$$
	z\to \tilde{f}(z):=f^c(z)I_{\{z\notin \R\setminus(\underline x,x)\}}+\tilde{l}(z)I_{\{z\in[\underline x,x]\}}
	$$
	is a convex minorant of $f$, and therefore $f^c\geq\tilde f$ on $\R$. But for $z\in I\cap(\underline x,x)$ we have that
	$$
	\tilde{f}(z)=\tilde{l}(z)>l^{f^c}_x(z)=f^c(z),
	$$
	a contradiction.
	
	Now suppose that $f^c<\tilde{l}$ on $(-\infty,x)$. Then
	$$
	z\to \overline{f}(z):=f^c(z)I_{\{z>x\}}+\tilde{l}(z)I_{\{z\leq x\}}
	$$
	is a convex minorant of $f$, and thus $f^c\geq\overline{f}$ on $\R$. But again, for $z\in I\cap(\infty,x)$ we have that
	$$
	\overline{f}(z)=\tilde{l}(z)>l^{f^c}_x(z)=f^c(z),
	$$
	which gives a required contradiction.
	
	3. Suppose $g^c(x)<f^c(x)$. The proof uses the arguments of the previous case (reverse the roles of $f^c$ and $g^c$, and consider $l^{f^c}_x,l^{g^c}_x$ on $[x,\infty)$).
\end{proof}

\begin{lemma}\label{lem:converge_hull}
Consider a sequence of measurable functions $f_n:\R\to\R$, $n\geq 1$. Suppose $f_n\downarrow f$ pointwise, for some measurable $f:\R\to\R$. Then $f^c_n\downarrow f^c$ pointwise as $n\to\infty$.  
\end{lemma}
\begin{proof}
	Fix $k\in\R$. Since $(f_n(k))_{n\geq1}$ is decreasing and bounded by $f(k)$, $\lim_{n\to\infty}f_n(k)$ exists. The same is true for the corresponding convex hulls. In particular, $\lim_{n\to\infty}f^c_n(k)\geq f^c(k)$. On the other hand,
	$$
	\lambda f(a)+(1-\lambda)f(b)=\lim_{n\to\infty}[\lambda f_n(a)+(1-\lambda)f_n(b)]\geq\lim_{n\to\infty}f^c_n(k),
	$$
	for all $a,b\in\R$ with $a\leq k\leq b$ and $\lambda\in[0,1]$ such that $\lambda a+(1-\lambda)b=k$. Taking infimum over all such $a,b,\lambda$ we obtain $f^c(k)\geq\lim_{n\to\infty}f^c_n(k)$.
\end{proof}

\section{Proofs of Section \ref{sec:mainExistence}}\label{sec:AppProofs}
\begin{proof}[Proof of Lemma \ref{lem:continuousProjection}]
	We first prove part (2). For a (Borel) measurable function $f:\R\to\R$, define $f^T:\R\to\R$ by
	$$
	f^T(x)=\int_\R f(z)dD_T(x,z),\quad x\in\R.
	$$
	Note that $\int_\R f(x) d\mu_n D_T(x)=\int_\R f^T(x)d\mu_n(x)$, and similarly for $\mu$ and $\mu D_T$. Hence to conclude that $\mu_nD_T\to\mu D_T$ it is enough to show that $\mu_nD_T\xrightarrow{w}\mu D_T$ and
	$$
	\int_\R f^T_0d\mu_n\to	\int_\R f^T_0d\mu,
	$$
	where, for each $t\in\R$, $f_t(x)=\lvert t-x\lvert$, $x\in\R$.
	
	Represent $\R\setminus T$ as a union of disjoint open intervals $\bigcup_{k\geq1}(l_k,r_k)$. Note that $l_k,r_k\in T$ for all $k\geq1$.
		
	We first establish the weak convergence. Let $f:\R\to\R$ be continuous and bounded. If $x\in T$, then $f^T(x)=f(x)$. If $x\notin T$, then $x\in(l_k,r_k)$ for some $k\geq1$, and then
	$$
	f^T(x)=\frac{r_k-x}{r_k-l_k}f(l_k)+\frac{x-l_k}{r_k-l_k}f(r_k)=L^f_{l_k,r_k}(x).
	$$
It follows that $f^T$ is also continuous and bounded, and therefore $\int_\R f^Td\mu_n\to\int_\R f^Td\mu$, which establishes the weak convergence.
	
	Now we deal with the convergence of first moments. Since $\mu_n\to\mu$ (i.e., w.r.t. $\sT_1$), it is enough to show that $f^T_0$ is continuous with at most linear growth.

First suppose that $\sup T=\infty$. Note that, if $0\in T$, then $f_0=f^T_0$ and we are done. On the other hand, if $0\notin T$, then $0\in(l_k,r_k)$ for some $k\geq 1$. We have that $f_0$ is linear on $(-\infty,l_k)\cup(r_k,\infty)$ and therefore $f^T_0=f_0$ on $(-\infty,l_k)\cup(r_k,\infty)$. It follows that $f^T_0=\max\{f_0,L^{f_0}_{l_k,r_k}\}$ and thus $f^T_0$ remains continuous with at most linear growth.

Now suppose that $\sup T<\infty$. If $\sup T\leq 0$, then $f^T_0(x)=(\sup T-x)^+-\sup T$, $x\in\R$. If $0<\sup T$ and $0\in T$, then $f^T_0=f_0$ on $(-\infty,\sup T]$ and $f^T_0\equiv\sup T$ on $(\sup T,\infty)$. Finally suppose that $0<\sup T$ and $0\notin T$. Then $0\in(l_k,r_k)$ for some $k\geq 1$ and $r_k\leq \sup T$. It follows that $f^T_0=f_0$ on $(-\infty,l_k]\cup[r_k,\sup T]$, $f^T_0=L^{f_0}_{l_k,r_k}$ on $(l_k,r_k)$ and $f^T_0\equiv\sup T$ on $(\sup T,\infty)$. It is evident that in all the cases $f^T_0$ remains continuous with at most linear growth.
	
	We now prove the uniqueness part (i.e., part (1)) .
	
	We first recall the irreducible decomposition of two measures $\mu\leq_{cd}\nu$, see Lemma \ref{lem:irreducible}. Let $x^*:=\sup\{k\in\R:P_\mu(k)=P_\nu(k)\}\in[-\infty,+\infty]$ with convention $\inf\emptyset=-\infty$. Represent an open set $\{k\in\R:P_\mu(k)<P_\nu(k)\}\cap(-\infty,x^*)$ by $\bigcup_{k\geq0} I_k=\bigcup_{k\geq0}(a_k,b_k)$, where $I_0=(x^*,\infty)$, and set $I_{-1}=\R\setminus\bigcup_{k\geq0} I_k$. If $\mu_k=\mu\lvert_{I_k}$, then there exists a unique decomposition $\nu=\sum_{k\geq-1}\nu_k$ such that
	$$
	\mu_{-1}=\nu_{-1},\quad\mu_0\leq_{cd}\nu_0\quad\textrm{and}\quad\mu_k\leq_c\nu_k\textrm{ for all } k\geq1. 
	$$
	In particular, any $\pi\in\Pi_S(\mu,\nu)$ admits a unique decomposition $\pi=\sum_{k\geq-1}\pi_k$ such that $\pi_0\in\Pi_S(\mu_0,\nu_0)$ and $\pi_k\in\Pi_M(\mu_k,\nu_k)$ for all $k\neq0$.
	
	Now let $\nu=\mu D_T$. Note that $P_{\mu D_T}(k)=P_\mu(k)$ for all $k\in T$. It follows that $\sup T\leq x^*$.
	
	By applying the arguments of Beiglb\"{o}ck and Juillet \cite[Proposition 4.1]{BeiglbockJuillet:16} to each $\pi_k\in\Pi_M(\mu_k,\nu_k)$ (for all $k\neq 0$), we obtain $\pi_k=\pi_{\mu_k,T}$. 
	
	We are left to show that $\pi_0=\pi_{\mu_0,T}$. Note that, if $x^*=\infty$, then we must have that $P_\mu$ and $P_\nu$ have the same asymptotic behaviour at $\infty$. In particular, $\overline\nu=\overline\mu$ and thus $\mu\leq_{c}\nu$. In this case the proof of $\pi_0=\pi_{\mu_0,T}$ is covered by the previous paragraph.
	
	If $\sup T=\infty$, then $x^*=\sup T=\infty$, and there is nothing to prove.
	
	Suppose $\sup T<\infty.$ Recall that $\sup T\leq x^*$. Since $\mu D_T((\sup T,\infty))=0$ we can, without loss of generality, assume that $\sup T=x^*$. (Indeed, if $\sup T<x^*$, then since $P_{\mu D_T}$ has slope $\mu(\R)$ to the right of $\sup T$ and $P_{\mu D_T}\geq P_\mu$ everywhere, we must have that $P_{\mu D_T}=P_\mu$ on $[\sup T,x^*]$. But then $x^*=\infty$, and again there is nothing to prove.) It follows that $\nu_0$ is an atomic measure concentrated on $x^*$, and therefore we necessarily have that $\pi_0=\pi_{\mu_0,T}$.
\end{proof}
\bibliographystyle{plainnat}

\end{document}